\documentclass[number,citesort,MSNbibl,seceqn,dvips]{arxbj}
\usepackage{mathrsfs}
\aid{0}
\volume{18}
\issue{1}
\pubyear{2012}
\firstpage{252}
\lastpage{278}
\doi{10.3150/10-BEJ335}

\makeatletter

\let\epsilon\varepsilon
\newcommand{\eqref}[1]{(\ref{#1})}

\newtheorem{theo}{Theorem}[section]
\newtheorem{prop}[theo]{Proposition}
\newtheorem{lem}[theo]{Lemma}
\newtheorem{cor}[theo]{Corollary}
\newremark{rem}[theo]{Remark}
\newproclaim{defi}[theo]{Definition}

\newcommand{\onrarrow}[1]{\mathrel{\vbox{\m@th\ialign{##\crcr
  $\hfil\scriptstyle\ #1\ \ \hfil$\crcr\noalign{\kern0.5pt\nointerlineskip}%
  \rightarrowfill\crcr}}}}

\makeatother

\begin{document}
\begin{frontmatter}

\title{On adaptive resampling strategies for sequential Monte Carlo methods}
\runtitle{Adaptive SMC}

\begin{aug}
\author[a]{\fnms{Pierre} \snm{Del Moral}\corref{}\thanksref{a}\ead[label=e1]{pierre.del-moral@inria.fr}},
\author[b]{\fnms{Arnaud} \snm{Doucet}\thanksref{b}\ead[label=e2]{a.doucet@stat.ubc.ca}} \and
\author[c]{\fnms{Ajay} \snm{Jasra}\thanksref{c}\ead[label=e3]{staja@nus.edu.sg}}
\runauthor{P. Del Moral, A. Doucet and A. Jasra}
\address[a]{Centre INRIA Bordeaux et Sud-Ouest \& Institut de Math\'{e}matiques de Bordeaux, Universit\'{e} de Bordeaux~I, 33405, France.
\printead{e1}}
\address[b]{Department of Statistics, University of British Columbia,
Vancouver BC, Canada V6T 1Z2.\\ \printead{e2}}
\address[c]{Department of Statistics and Applied Probability, National
University of Singapore, Singapore 117546, Singapore. \printead{e3}}
\end{aug}

\received{\smonth{5} \syear{2009}}
\revised{\smonth{9} \syear{2010}}

%
\begin{abstract}
Sequential Monte Carlo (SMC) methods are a class of techniques to sample
approximately from any sequence of probability distributions using a
combination of importance sampling and resampling steps. This paper is
concerned with the convergence analysis of a class of SMC methods where the
times at which resampling occurs are computed online using criteria
such as
the effective sample size. This is a popular approach amongst practitioners
but there are very few convergence results available for these methods. By
combining semigroup techniques with an original coupling argument, we
obtain functional central limit theorems and uniform exponential
concentration estimates for these algorithms.
\end{abstract}

%

\end{frontmatter}

\section{Introduction}
\label{sec:intro}

Sequential Monte Carlo (SMC) methods are a generic class of simulation-based
algorithms to sample approximately from any sequence of probability
distributions. These methods are now extensively used in engineering,
statistics and physics; see \cite{cappe2005,fk,arnaud,liu2001} for many
applications. Sequential Monte Carlo\ methods approximate the target
probability distributions
of interest by a large number of random samples, termed particles, which
evolve over time according to a combination of importance sampling and
resampling steps.

In the resampling steps, new particles are sampled with replacement
from a
weighted empirical measure associated to the current particles; see
Section \ref{adaptsec} for more details. These resampling steps are crucial and,
without them, it is impossible to obtain time uniform convergence results
for SMC\ estimates. However, resampling too often has a~negative effect as
it decreases the number of distinct particles. Hence, a resampling
step
should only be applied when necessary. Consequently, in most practical
implementations of SMC, the times at which resampling occurs are
selected by
monitoring a criterion that assesses the quality of the current particle
approximation. Whenever this criterion is above or below a given
threshold, a resampling step is triggered. This approach was originally proposed
in \cite{liuchen95} and has been widely adopted ever since \cite
{cappe2005}, Section~7.3.2.

For this class of adaptive SMC\ methods, the resampling times are computed
online using our current SMC approximation and thus are random. However,
most of the theoretical results on SMC algorithms assume resampling occurs
at deterministic times; see~\cite{douc} for an exception discussed later.
The objective of this paper is to provide convergence results for this
type of adaptive SMC\ algorithm. This is achieved using a coupling
argument. Under some assumptions, the random resampling times converge
almost surely as the number of particles goes to infinity toward some
deterministic (but not explicitly known) resampling times. We show here that
the difference, in probability, between the reference SMC algorithm
based on
these deterministic but unknown resampling times and the adaptive SMC\
algorithm is exponentially small in the number of particles. This
allows us
to straightforwardly transfer the convergence results of the reference
SMC\
algorithm to the adaptive SMC\ algorithm. In particular, we establish
functional central limit theorems and new exponential concentration
estimates that improve over those presented in \cite{fk}, Section 7.4.3.
Note that some exponential concentration estimates have also been
established in \cite{cappe2005}, Theorem 9.4.12, using different techniques
and a weaker assumption. The constants appearing in \cite%
{cappe2005}, Theorem 9.4.12, are not explicit so the comparison between
these two results is
difficult. In a specific example, we found our bound to be significantly
tighter but have not established it in a general case.

The rest of the paper is organized as follows. In Section \ref%
{sec:mainresult}, we present the class of adaptive SMC\ algorithms studied
here and our main coupling result. A\ precise description of the
sequence of
distributions approximated by the reference SMC\ algorithm is given in
Section~\ref{dmosec} and a theoretical analysis of the reference SMC\
algorithm is presented in Section~\ref{adpative}. In particular, we propose
an original concentration analysis to obtain exponential estimates for
SMC\
approximations. These results are used to obtain a concentration result for
the empirical criteria around their limiting values. The results above are
used, in Section \ref{asymptsec}, to bound the differences between the
deterministic resampling times and their empirical approximations, up
to an
event with an exponentially small probability. Finally, we analyze the
fluctuations of adaptive SMC algorithms in Section \ref{functtclsec}.

\section{Adaptive SMC\ algorithms and main results}
\label{sec:mainresult}

\subsection{Notation and conventions}

Let $\mathcal{M}(E)$, $\mathcal{P}(E)$ and $\mathcal{B}_{b}(E)$ denote,
respectively, the set of bounded and signed measures, the subset of all
probability measures on some measurable space $(E,\mathscr{E})$ and the
Banach space of all bounded and measurable functions $f$ on $E$ when
equipped with norm $\Vert f\Vert=\sup_{x\in E}{|f(x)|}$.
$\operatorname{Osc}_{1}(E)$
is the set of $\mathscr{E}$-measurable functions $f$ with oscillations
$%
\operatorname{osc}(f)=\sup_{(x,y)\in E^{2}}{\{|f(x)-f(y)|\}}\leq1$.
$\mu(f)=\int
\mu(\mathrm{d}x)f(x)$ is the integral of a~function $f\in\mathcal{B}_{b}(E)$,
w.r.t. a measure $\mu\in\mathcal{M}(E)$. $\mu(A)=\mu(1_{A})$ with
$A\in%
\mathscr{E}$ and $1_{A}$ the indicator of $A$. $\delta_{a}$ is the Dirac
measure. A bounded integral operator $M$ from a~measurable space $(E,%
\mathscr{E})$ into another $(F,\mathscr{F})$ is an operator $f\mapsto M(f)$
from $\mathcal{B}_{b}(F)$ into~$\mathcal{B}_{b}(E)$ such that the
functions $%
M(f)(x)=\int_{F}M(x,\mathrm{d}y)f(y)$ are measurable and bounded for any $f\in
\mathcal{B}_{b}(F)$. A bounded integral operator $M$ from $(E,\mathscr{E})$
into $(F,\mathscr{F})$ also generates a~dual operator $\mu\mapsto\mu M$
from $\mathcal{M}(E)$ into $\mathcal{M}(F)$ defined by $(\mu
M)(f):=\mu
(M(f))$. If constants are written with an argument, then they depend
only on
this given argument. The tensor product of functions is written
$\otimes$.
For any generic sequence $ \{ z_{n} \} _{n\geq0}$, we denote $%
z_{i:j}= ( z_{i},z_{i+1},\ldots,z_{j} ) $ for $i\leq j$.

\subsection{Adaptive sequential Monte Carlo methods}

\label{adaptsec}

SMC methods are a popular class of methods for sampling random variables
distributed approximately according to the Feynman--Kac path measures
%
\begin{eqnarray}\label{form1}
\eta^*_{n}(f_{n}) &=&\gamma^*_{n}(f_{n})/\gamma^* _{n}(1)
\qquad \mbox{with } \gamma^*_{n}(f_{n})=\mathbb{E} (f_{n} ( X_{0:n} )
W_{0:n-1} ( X_{1:n-1} ) ) , \\
\label{form2}
\widehat{\eta}^*_{n}(f_{n}) &=&\widehat{\gamma
}^*_{n}(f_{n})/\widehat{\gamma}%
^*_{n}(1) \qquad \mbox{with } \widehat{\gamma}^*_{n}(f_{n})=\mathbb{E}%
(f_{n} ( X_{0:n} ) W_{0:n} ( X_{1:n} ) ) ,
\end{eqnarray}
where $(X_{n})_{n\geq0}$ is a Markov chain on $(E_{n},\mathscr{E}%
_{n})_{n\geq0}$ with transition kernels $ ( M_{n} ) _{n>0}$, $%
(G_{n})_{n>0}$ is a sequence of non-negative potential functions on $ (
E_{n} ) _{n>0}\ $and the importance weight function is defined by
%
\begin{equation}\label{eq:weight}
W_{p,q}\dvtx x_{p+1:q}\in E_{p+1}\times\cdots\times E_{q}\mapsto
W_{p,q}(x_{p+1:q}):=\prod_{p<k\leq q}G_{k}(x_{k}).
\end{equation}

The basic SMC\ method proceeds as follows. Given $N$ particles distributed
approximately according to $\eta^*_{n-1}$, these particles first evolve
according to the transition kernel~$M_{n}$. In a second stage, particles
with low relative $G_{n}$-potential value are killed and those with a larger
relative potential are duplicated. However, as noted in Section \ref%
{sec:intro}, resampling at each time step is wasteful and should only be
performed when necessary.

This has motivated researchers to introduce new resampling strategies where
the resampling step is only triggered when a criterion is satisfied;
this is
typically computed via the current particle approximation (see Section
\ref%
{sec:adaptivecriteria}). Such adaptive SMC\ algorithms proceed as follows.
Let $t_{n}^{N}$ denote the $n$th resampling time of the adaptive
SMC\ algorithm. After the $n$th resampling step, assume we have
the following empirical measure approximation of $\widehat{\eta}%
^{*N}_{t_{n}} $ denoted
\[
\widehat{\eta}_{t_{n}^{N}}^{*N} ( \cdot ) :=\frac{1}{N}%
\sum_{i=1}^{N}\delta_{\widehat{\mathcal{Y}}_{n}^{(N,i)}}(\cdot),
\]
where $\widehat{\mathcal{Y}}_{n}^{(N,i)}:=\widehat
{Y}_{0:t_{n}^{N}}^{ (
N,i ) }$. We propagate forward these $N$ paths by generating $%
Y_{t_{n}^{N}+1:t_{n+1}^{N}}^{(N,i)}$ according to the transition kernel
$%
M_{t_{n}^{N}+1:t_{n+1}^{N}}:=M_{t_{n}^{N}+1}M_{t_{n}^{N}+2}\cdots
M_{t_{n+1}^{N}}$ of the reference Markov chain initialized at $\widehat
{Y}%
_{t_{n}^{N}}^{(N,i)}$, up to the first time ($t_{n+1}^{N}$) the importance
weights of the $N$ path samples given by $%
W_{t_{n}^{N},t_{n+1}^{N}}(Y_{t_{n}^{N}+1:t_{n+1}^{N}}^{(N,i)})$ become, in
some sense, degenerate.

At time $t_{n+1}^{N}$ the weighted occupation measure of the system
\[
\widetilde{\eta}_{t_{n+1}^{N}}^{*N}(\cdot):=\sum_{i=1}^{N}\frac{%
W_{t_{n}^{N},t_{n+1}^{N}}(Y_{t_{n}^{N}+1:t_{n+1}^{N}}^{(N,i)})}{%
\sum
_{j=1}^{N}W_{t_{n}^{N},t_{n+1}^{N}}(Y_{t_{n}^{N}+1:t_{n+1}^{N}}^{(N,j)})}
\delta_{\mathcal{Y}_{n+1}^{(N,i)}}(\cdot)
\]
is a particle approximation of $\widehat{\eta}^*_{t_{n+1}^{N}}(\cdot),$
where $\mathcal{Y}_{n+1}^{(N,i)}:=(\widehat{\mathcal{Y}}%
_{n}^{(N,i)},Y_{t_{n}^{N}+1:t_{n+1}^{N}}^{(N,i)})$. After the resampling
step, this measure is replaced by an empirical measure
\[
\widehat{\eta}_{t_{n+1}^{N}}^{*N} ( \cdot ) :=\frac{1}{N}%
\sum_{i=1}^{N}\delta_{\widehat{\mathcal{Y}}_{n+1}^{(N,i)}}(\cdot)
\]
associated with $N$ path particles $\widehat{\mathcal
{Y}}_{n+1}^{(N,i)}:=%
\widehat{Y}_{0:t_{n+1}^{N}}^{(N,i)}$ that are resampled from
$\widetilde{\eta%
}_{t_{n+1}^{N}}^{*N}$; see, for example, \cite{arnaud} for alternative
resampling
schemes.

\vspace*{-2pt}\subsection{Some empirical criteria}\vspace*{-2pt}
\label{sec:adaptivecriteria}

Two well-known criteria used in the SMC\ literature to trigger the
resampling mechanism are now discussed. In both cases, the resampling
times $%
( t_{n}^{N} ) _{n\geq0}$ are random variables that depend on the
current SMC\ approximation.

$\bullet$ \textit{Squared coefficient of variation}. After the resampling
step at time $t_{n}^{N}$, the particles explore the state space up to the
first time $(s=t_{n+1}^{N}$) the squared coefficient of variation of the
unnormalized weights is larger than some prescribed threshold $a_{n}$
%
\begin{equation} \label{eq:empiricalvariance}
C_{t_{n}^{N},s}^{N}=\frac{1}{N}\sum_{i=1}^{N} \Biggl(
W_{t_{n}^{N},s}\bigl(Y_{t_{n}^{N}+1:s}^{(N,i)}\bigr)\Big/\frac{1}{N}%
\sum_{j=1}^{N}W_{t_{n}^{N},s}\bigl(Y_{t_{n}^{N}+1:s}^{(N,j)}\bigr) \Biggr) ^{2}-1\geq
a_{n}.
\end{equation}
This is equivalent to resampling when the effective sample size (ESS),
defined as $\mathit{ESS}=N ( 1+C_{t_{n}^{N},s}^{N} ) ^{-1}$, is below a
prescribed threshold as proposed in \cite{liuchen95}.

$\bullet$ \textit{Entropy}. After the resampling step at time $t_{n}^{N}$,
the particles explore the state space up to the first time $(s=t_{n+1}^{N}$)
the relative entropy of the empirical particle measure w.r.t. its weighted
version is larger than some threshold $a_{n}$
%
\begin{equation}\label{eq:empiricalentropy}
C_{t_{n}^{N},s}^{N}:=-\frac{1}{N}\sum_{i=1}^{N}\log
W_{t_{n}^{N},s}\bigl(Y_{t_{n}^{N}+1:s}^{(N,i)}\bigr)\geq a_{n}.
\end{equation}
\subsection{Statement of some results}
\label{sec:statement}

The following section provides a guide of the major definitions and results
in this paper; these will be repeated at the relevant stages in the\vadjust{\goodbreak}
paper.

\subsubsection{A limiting reference SMC\ algorithm}\vspace*{-2pt}
\label{sec:limiting}

Let $(t_{n})_{n\geq0}$ be the deterministic sequence of time steps obtained
by replacing the empirical criteria $C_{t_{n}^{N},s}^{N}$ by their limiting
values $C_{t_{n},s}$ as $N\uparrow\infty$, that is, $t_{n+1}:=\inf{ \{
t_{n}<s\dvt C_{t_{n},s}\geq a_{n} \} }.$ In all situations, set $%
t_{0}^{N}=t_{0}=0$. For the criterion (\ref{eq:empiricalvariance}), the
limiting criterion $C_{t_{n},s}$ is given by
%
\begin{equation}\label{eq:limitingvariance}
\frac{\mathbb{E}_{t_{n},\widehat{\eta}_{t_{n}}^{\ast}}(W_{t_{n},s} (
X_{t_{n}+1:s} ) ^{2})}{\mathbb{E}_{t_{n},\widehat{\eta
}_{t_{n}}^{\ast
}} ( W_{t_{n},s} ( X_{t_{n}+1:s} ) ) ^{2}}-1,
\end{equation}
whereas for (\ref{eq:empiricalentropy}) it is given by
%
\begin{equation}\label{eq:limitingentropy}
-\mathbb{E}_{t_{n},\widehat{\eta}_{t_{n}}^{\ast}} (\log
W_{t_{n},s} ( X_{t_{n}+1:s} ) ).
\end{equation}
Here, $\mathbb{E}_{t_{n},\widehat{\eta}_{t_{n}}^{\ast}}$ is the
expectation w.r.t. the law $\mathbb{P}_{t_{n},\widehat{\eta
}_{t_{n}}^{\ast}}
$ of the random path of variables that starts at time $t_{n}$ at the end
point $X_{t_{n}}=\widehat{X}_{t_{n}}$ of $\widehat{\mathcal{X}}_{n}:=
\widehat{X}_{0:t_{n}}$ distributed according to~$\widehat{\eta}%
_{t_{n}}^{\ast}$ and evolves according to the Markov kernels
$M_{t_{n}+1:s}$. In \cite{cornebise2008}, the limiting expression for the normalized
effective sample size $N^{-1}\mathit{ESS}$ has been established. An alternative
entropy criterion has also been proposed that can be applied when the
potential functions~$(G_{n})_{n>0}$ are not strictly positive on $ (
E_{n} ) _{n>0}$.

\subsubsection{An exponential coupling theorem}\vspace*{-2pt}

We give our main results, which hold under the following regularity
condition:
%
\begin{equation}\label{condGprime}
(G)\quad \forall n\geq1\qquad  q_{n}^{\prime}:=\sup_{(x,y)\in
(E_{n})^{2}} \bigl( {G_{n}(x)}/{G_{n}(y)} \bigr)<\infty.
\end{equation}
We refer the reader to \cite{fk}, Chapter 3, for a thorough discussion
in the
case where $(G)$ does not apply. To state our results, we first require the
following definition.\vspace*{-2pt}

\begin{defi}
Let $\mathcal{Y}_{n}^{(N)}:=(\mathcal{Y}_{n}^{(N,1)},\mathcal{Y}%
_{n}^{(N,2)},\ldots,\mathcal{Y}_{n}^{(N,N)})$ and $\widehat{\mathcal
{Y}}%
_{n}^{ ( N ) }:=(\widehat{\mathcal{Y}}_{n}^{(N,1)},\widehat{%
\mathcal{Y}}_{n}^{(N,2)},\ldots,\allowbreak\widehat{\mathcal{Y}}_{n}^{(N,N)})$ denote
the $N$ particles associated to the adaptive SMC\ algorithm resampling at
times $ ( t_{n}^{N} ) _{n\geq0}$ and let $\mathcal{X}_{n}^{(N)}:=(%
\mathcal{X}_{n}^{(N,1)},\mathcal{X}_{n}^{(N,2)},\ldots,\mathcal{X}%
_{n}^{(N,N)})$ and $\widehat{\mathcal{X}}_{n}^{ ( N ) }:=(\widehat{%
\mathcal{X}}_{n}^{(N,1)},\widehat{\mathcal{X}}_{n}^{(N,2)},\allowbreak\ldots
,\widehat{%
\mathcal{X}}_{n}^{(N,N)})$ denote the $N$ particles associated to the
reference SMC\ algorithm resampling at times $ ( t_{n} ) _{n\geq0}$%
. We also suppose that $(\mathcal{X}_{n}^{ ( N ) },\widehat{%
\mathcal{X}}_{n}^{(N)})$ and $(\mathcal{Y}_{n}^{ ( N ) },\widehat{%
\mathcal{Y}}_{n}^{(N)})$ coincide on every time interval $0\leq n\leq m$,
once $t_{n}^{N}=t_{n}$, for every $0\leq n\leq m$. This condition
corresponds to the coupling of the two processes on the event $\bigcap
_{0\leq
n\leq m}\{t_{n}^{N}=t_{n}\}$.\vspace*{-2pt}
\end{defi}

The first result is a non-asymptotic exponential concentration estimate.
The probability measures $\eta_{n}^{N}$ and $\eta_{n}$ are
introduced below. They can be thought of as analogues of~$\eta
_{t_{n}}^{\ast N}$
and $\eta_{t_{n}}^{\ast}$; see Section \ref{dmosec} for formal definitions.\vspace*{-2pt}

\begin{theo}
\label{theo2} For any $n\geq0$, $f_{n}\in\operatorname{Osc}_{1}(E_{0}\times
\cdots\times E_{t_{n}})$, any $N\geq1$ and any $0\leq\epsilon\leq1/2$,
there exist $c_{1}<\infty$, $0<c_{2}(n)<\infty$ such that we have the
exponential concentration estimate
\[
\mathbb{P} \bigl( \vert [ \eta_{_{n}}^{N}-\eta_{_{n}} ]
(f_{n}) \vert\geq\epsilon \bigr) \leq c_{1}\exp\{-{N\epsilon^{2}%
}/{c_{2}(n)}\}\vadjust{\goodbreak}
\]
for the empirical measures $\eta_{_{n}}^{N}(\cdot)=\frac{1}{N}%
\sum_{i=1}^{N}\delta_{ ( \widehat{\mathcal{X}}%
_{n-1}^{(N,i)},X_{t_{n-1}+1:t_{n}}^{(N,i)} ) }(\cdot).$ In addition,
under appropriate regularity conditions on $ ( M_{k} ) _{k>0}$ and $%
( G_{k} ) _{k>0}$ given in Section \ref{introconcc} the above
estimates are valid for the marginal measures associated to the time
parameters $t_{n-1}\leq p\leq t_{n}$ for some constant $c_{2}(n)=c_{2}$.
\end{theo}

The second result is an exponential coupling theorem.

\begin{theo}
\label{theo1} Assume the threshold parameters $(a_{n})_{n\geq0}$ are
sampled realizations of a collection of absolutely continuous random
variables $A=(A_{n})_{n\geq0}$. Then, for almost every realization of the
sequence $(A_{n})_{n\geq0}$, $(\mathcal{X}_{n}^{(N)},\widehat
{\mathcal{X}}%
_{n}^{ ( N ) })_{n\geq0}$ and $(\mathcal{Y}_{n}^{(N)},\widehat{%
\mathcal{Y}}_{n}^{ ( N ) })_{n\geq0}$ are such that, for every $%
m\geq0$ and any $N\geq1$, there exist $0<c_{1}(m),c_{2}(m)<\infty$ and
almost surely $\epsilon(m,A)\equiv\epsilon(m)>0$ such that
\[
\mathbb{P} \bigl( \exists0\leq n\leq m\ \bigl(\mathcal{Y}_{n}^{(N)},\widehat{%
\mathcal{Y}}_{n}^{ ( N ) }\bigr)\not=\bigl(\mathcal{X}_{n}^{(N)},\widehat{%
\mathcal{X}}_{n}^{ ( N ) }\bigr)\vert A \bigr) \leq c_{1}(m)\mathrm{e}^{-N\epsilon
^{2}(m)/c_{2}(m)}.
\]
\end{theo}

Up to an event having an exponentially small occurrence probability,
Theorem %
\ref{theo1} allows us to transfer many estimates of the reference SMC
algorithm $(\mathcal{X}_{n}^{ ( N ) },\widehat{\mathcal{X}}%
_{n}^{ ( N ) })_{n\geq0}$ resampling at deterministic times to the
adaptive SMC\ algorithm $(\mathcal{Y}_{n}^{(N)},\widehat{\mathcal{Y}}_{n}^{ ( N ) })_{n\geq0}$.

The proofs of Theorems~\ref{theo2} and~\ref{theo1} are detailed,
respectively, in Sections \ref{conctheo} and~\ref{randomizedr}.

\section{Description of the models}
\label{dmosec}


\subsection{Feynman--Kac distributions flow}
\label{secgen}

We consider a sequence of measurable state spaces $ ( S_{n},\mathscr{S}
_n ) _{n\geq0}$, a probability measure $\eta_{0}\in\mathcal{P}%
(S_{0}) $ and a sequence of Markov transitions $\mathcal{M}%
_{n}(x_{n-1},\mathrm{d}x_{n})$ from $S_{n-1}$ into $S_{n}$ for $n\geq1$. Let $ (
\mathcal{X}_{n} )_{n\geq0} $ be a Markov chain with initial distribution
$\operatorname{Law}(\mathcal{X}_{0})=\eta_{0}$ and elementary
transitions $%
\mathbb{P}(\mathcal{X}_{n}\in \mathrm{d}y|\mathcal{X}_{n-1}=x)=\mathcal{M}%
_{n}(x,\mathrm{d}y).$ Let $ (\mathcal{G}_{n} )_{n\geq0}$ be a sequence of
non-negative and bounded potential functions on $S_{n}$. To simplify the
presentation, and to avoid unnecessary technicalities, it is supposed $%
\mathcal{G}_{n}\in(0,1)$ for $n\geq1$ with
%
\begin{equation}\label{condG}
(\mathcal{G})\quad  q_{n}:=\sup_{(x,y)\in S_{n}^{2}} \bigl( {\mathcal{G}%
_{n}(x)}/{\mathcal{G}_{n}(y)} \bigr) <\infty.
\end{equation}
The Boltzmann--Gibbs transformation $\Psi_{n}$ associated to $\mathcal
{G}%
_{n} $ is the mapping
\[
\Psi_{n}\dvtx\mu\in\mathcal{P}(S_{n})\mapsto\Psi_{n}(\mu)\in
\mathcal{P}%
(S_{n}) \qquad \mbox{with } \Psi_{n}(\mu)(\mathrm{d}x):=\frac{1}{\mu(%
\mathcal{G}_{n})}\mathcal{G}_{n}(x)\mu(\mathrm{d}x).
\]
Notice that $\Psi_{n}(\mu)$ can be rewritten as a nonlinear Markov
transport equation
\[
\Psi_{n}(\mu)(\mathrm{d}y)= ( \mu\mathcal{S}_{n,\mu} )
(\mathrm{d}y)=\int_{S_{n}}\mu(\mathrm{d}x)\mathcal{S}_{n,\mu}(x,\mathrm{d}y)
\]
with\vadjust{\goodbreak} $\mathcal{S}_{n,\mu}(x,\mathrm{d}y):=\mathcal{G}_{n}(x)\delta
_{x}(\mathrm{d}y)+(1-\mathcal{G}_{n}(x))\mathcal{G}_{n}(y)\mu(\mathrm{d}y)/\mu(\mathcal{G}_{n}).$

Let $(\eta_{n},\widehat{\eta}_{n})_{n\geq0}$ be the flow of probability
measures, both starting at $\eta_{0}=\widehat{\eta}_{0}$, and
defined for
any $n\geq1$ by the following recursion
%
\begin{equation}\label{predeq}
\forall n\geq0\qquad  \eta_{n+1}=\widehat{\eta}_{n}\mathcal{M}_{n+1} %
\qquad \mbox{with } \widehat{\eta}_{n}:=\Psi_{n}(\eta_{n})=\eta_{n}%
\mathcal{S}_{n,\eta_{n}}.
\end{equation}
It can be checked that the solution $(\eta_{n},\widehat{\eta}_{n})$ of
these recursive updating prediction equations have the following functional
representations:
%
\begin{equation}
\eta_{n}(f_{n})=\gamma_{n}(f_{n})/\gamma_{n}(1) \quad \mbox{and}\quad
\widehat{\eta}_{n}(f_{n})=\widehat{\gamma}_{n}(f_{n})/\widehat
{\gamma}%
_{n}(1) \label{qf}
\end{equation}
with the unnormalized Feynman--Kac measures $\gamma_{n}$ and $\widehat
{\gamma}_{n}$ defined by the formulae
%
\begin{equation}\label{unn}
\gamma_{n}(f_{n})=\mathbb{E} \biggl[ f_{n}(\mathcal{X}_{n})\prod_{0<k<n}%
\mathcal{G}_{k}(\mathcal{X}_{k}) \biggr]\quad  \mbox{and}\quad  \widehat{%
\gamma}_{n}(f_{n})=\gamma_{n}(f_{n}\mathcal{G}_{n}).
\end{equation}

\subsection{Feynman--Kac semigroups}
\label{secfk}

To analyze SMC\ methods, we introduce the Feynman--Kac
semigroup associated to the flow of measures $(\gamma_{n})_{n\geq0}$
and $%
(\eta_{n})_{n\geq0}$. Let us start by denoting by $\mathcal{Q}%
_{n+1}(x_{n},\mathrm{d}x_{n+1})$ the bounded integral operator from $S_{n}$ into
$S_{n+1}$ defined by
\[
\mathcal{Q}_{n+1}(x_{n},\mathrm{d}x_{n+1})=\mathcal{G}_{n}(x_{n})\mathcal{M}_{n+1}(x_{n},\mathrm{d}x_{n+1}).
\]
Let $(\mathcal{Q}_{p,n})_{0\leq p\leq n}$ be the corresponding linear
semigroup defined by $\mathcal{Q}_{p,n}=\mathcal{Q}_{p+1}\mathcal{Q}%
_{p+2}\cdots\mathcal{Q}_{n}$ with the convention $\mathcal
{Q}_{n,n}=I$, the
identity operator. Note that $\mathcal{Q}_{p,n}$ is alternatively
defined by
%
\begin{equation}\label{qpnref}
\mathcal{Q}_{p,n}(f_{n})(x_{p})=\mathbb{E} \biggl[f_{n}(\mathcal{X}%
_{n})\prod_{p\leq k<n}\mathcal{G}_{k}(\mathcal{X}_{k}) \big|\mathcal{X}%
_{p}=x_{p} \biggr].
\end{equation}
Using the Markov property, it follows that
\[
\gamma_{n}(f_{n})=\mathbb{E} \biggl[\mathbb{E} \biggl[f_{n}(\mathcal{X}%
_{n})\prod_{p\leq k<n}\mathcal{G}_{k}(\mathcal{X}_{k}) \big|\mathcal
{X}_{p}%
\biggr]\prod_{0<k<p}\mathcal{G}_{k}(\mathcal{X}_{k}) \biggr]=\gamma_{p} (
\mathcal{Q}_{p,n}(f_{n}) ) .
\]
The last assertion shows that $(\mathcal{Q}_{p,n})_{0\leq p\leq n}$ is the
semigroup associated with the unnormalized measures $(\gamma
_{n})_{n\geq
0} $. Denote its normalized version by%
%
\begin{equation}\label{qnormref}
\mathcal{P}_{p,n}(f_{n}):=\frac{\mathcal{Q}_{p,n}(f_{n})}{\mathcal
{Q}_{p,n}(1)}.
\end{equation}
Finally, denote by $(\Phi_{p,n})_{0\leq p\leq n}$ the nonlinear semigroup
associated to the flow of normalized measures $(\eta_{n})_{n\geq0}$:
$\Phi
_{p,n}=\Phi_{n}\circ\cdots\circ\Phi_{p+2}\circ\Phi_{p+1}$ with
the
convention \mbox{$\Phi_{n,n}=I$}, the identity\vadjust{\goodbreak} operator and $\Phi_{n}(\mu
)=\mu(%
\mathcal{G}_{n-1}\mathcal{M}_{n})/\mu(\mathcal{G}_{n-1})$, $\mu\in
\mathcal{P}(S_{n-1})$. Note that $(\Phi_{p,n})_{0\leq p\leq n}$ can be
alternatively defined in terms of $(\mathcal{Q}_{p,n})_{0\leq p\leq n}$
using
%
\begin{equation}
\Phi_{p,n}(\eta_{p})(f_{n})=\frac{\gamma_{p}\mathcal
{Q}_{p,n}(f_{n})}{%
\gamma_{p}\mathcal{Q}_{p,n}(1)}=\frac{\eta_{p}\mathcal
{Q}_{p,n}(f_{n})}{%
\eta_{p}\mathcal{Q}_{p,n}(1)}. \label{phiQ}
\end{equation}

\subsection{Path space and excursion models}
\label{secexc}

Let $ (X_{n} )_{n\geq0}$ be a Markov chain taking values in some
measurable state spaces $E_{n}$ with elementary transitions $%
M_{n}(x_{n-1},\mathrm{d}x_{n})$ and initial distribution $\eta_{0}=\operatorname{Law}
(X_{0})$. In addition, introduce a sequence of non-negative potential
functions $ ( G_{n} ) _{n>0}$ on the state spaces~$ (
E_{n} ) _{n>0}$. To simplify the presentation, it is assumed that $%
G_{n}\in ( 0,1 ) $ $\forall n>0$.

We associate to an increasing sequence of time parameters
$(t_{n})_{n\geq0}$
the excursion-valued random variables $X_{0}$ for $n=0$ and $%
X_{t_{n-1}+1:t_{n}}$ for $n\geq1$. We also define the random path
sequences%
\[
\mathcal{X}_{n}:=X_{0:t_{n}}\in E_{t_{n}}^{\prime}
\]
with the convention $E_{n}^{\prime}:=E_{0}\times\cdots\times E_{n}$. Note
that $ (\mathcal{X}_{n} )_{n\geq0}$ forms a Markov chain
%
\begin{equation}\label{refxa}
\mathcal{X}_{n+1}:= ( \mathcal{X}_{n},X_{t_{n}+1:t_{n+1}} )
\end{equation}
taking values in the excursion spaces $S_{n}:=E_{t_{n}}^{\prime}$. Now
adopting the potential functions
%
\begin{equation}\label{laformm}
\forall n\geq1 \qquad \mathcal{G}_{n}(\mathcal{X}_{n}):=W_{t_{n-1}:t_{n}}%
( X_{t_{n-1}+1:t_{n}} )
\end{equation}
in \eqref{unn}, we readily find that
\[
\gamma_{n}(f_{n})=\mathbb{E} \biggl[f_{n}(\mathcal{X}_{n})\prod_{0<k<n}%
\mathcal{G}_{k}(\mathcal{X}_{k}) \biggr]=\mathbb{E} [%
f_{n}(X_{0:t_{n}})W_{0:t_{n-1}} ( X_{1:t_{n-1}} ) ].
\]

By definition of the potential functions $\mathcal{G}_{n}$ of the excursion
Feynman--Kac model (\ref{laformm}), it is easily proved that the
condition $(%
\mathcal{G})$ (equation~(\ref{condG})) is satisfied as soon as $(G)$
introduced
in (\ref{condGprime}) holds true. More precisely, it holds that $(G)$
implies $(\mathcal{G})$ with
\[
q_{n}\leq M-\sup \biggl\{\frac{{%
W_{t_{n-1}:t_{n}} ( x_{t_{n-1}+1:t_{n}} ) }}{{W_{t_{n-1}:t_{n}}%
( y_{t_{n-1}+1:t_{n}} ) }} \biggr\}\qquad  \biggl(\leq\prod_{t_{n-1}<k\leq
t_{n}}q_{k}^{\prime} \biggr),
\]
where the essential supremum $M-\sup{ \{ \cdot\} }$ is
taken over all admissible paths $x_{t_{n-1}+1:t_{n}}$ and $%
y_{t_{n-1}+1:t_{n}}$ of the underlying Markov chain $(X_{n})_{n\geq0}$.

\subsection{Functional criteria}\label{adaptsc}

In Section \ref{secexc}, we have assumed that an increasing
sequence of time parameters $(t_{n})_{n\geq0}$ was available. We now
introduce the functional criteria used to build this sequence. To connect
the empirical criteria with their limiting functional versions, the latter
need to satisfy some weak regularity conditions that are given below.

\begin{defi}
We consider a sequence of functional criteria
\[
\forall n\geq0,\ \forall p\leq q \qquad \mathcal{H}_{p,q}^{(n)}\dvtx\mu
\in\mathcal{P} ( E_{q}^{\prime} ) \mapsto\mathcal{H}%
_{p,q}^{(n)} ( \mu ) \in\mathbb{R}_{+}
\]
satisfying the following Lipschitz type regularity condition
%
\begin{equation}\label{lipsH}
\bigl\vert\mathcal{H}_{p,q}^{(n)} ( \mu_{1} ) -\mathcal{H}%
_{p,q}^{(n)} ( \mu_{2} ) \bigr\vert\leq\delta
\bigl(H_{p,q}^{(n)}\bigr)\int| [ \mu_{1}-\mu_{2} ] (h)
|H_{p,q}^{(n)}(\mathrm{d}h)
\end{equation}
for some collection of bounded measures $H_{p,q}^{(n)}$ on $\mathcal{B}
_{b}(E_{q}^{\prime})$ such that
\[
\delta\bigl(H_{p,q}^{(n)}\bigr):=\int\operatorname
{osc}(h)H_{p,q}^{(n)}(\mathrm{d}h)<\infty.
\]
\end{defi}

We illustrate this construction with the pair of functional criteria
discussed in Section~\ref{sec:limiting}. When we consider (\ref%
{eq:limitingvariance}), the functional
%
\begin{equation}\label{eq:limitvariance}
\mathcal{H}_{p,q}^{(n)} ( \mu ) =\mu \biggl( \biggl[ \frac{W_{p,q}}{%
\mu(W_{p,q})}-1 \biggr] ^{2} \biggr)
\end{equation}
coincides with the squared coefficient of variation of the weights
w.r.t. $%
\mu$. When we consider (\ref{eq:limitingentropy}), the functional
%
\begin{equation}\label{eq:limitentropy}
\mathcal{H}_{p,q}^{(n)} ( \mu ) =\operatorname{Ent} ( \mathrm{d}\mu
|W_{p,q}\mathrm{d}\mu ) :=-\mu ( \log{W_{p,q}} )
\end{equation}
measures the relative entropy distance between $\mu$ and the updated
weighted measure. Under the condition $(\mathcal{G})$ stated in (\ref
{condG}%
), it is an elementary exercise to check that the above pair of criteria
satisfy (\ref{lipsH}). In the first case (\ref{eq:limitvariance}), we can
take \mbox{$H_{p,q}^{(n)}=c[\delta_{W_{p,q}^{2}}+\delta_{W_{p,q}}]$} for some
constant $c$ sufficiently large. In the second case (\ref{eq:limitentropy}),
we can take $H_{p,q}^{(n)}=c\delta_{W_{p,q}}$, again for some $c$ large
enough.

\subsection{Resampling times construction}

We now explain how to define the sequence of resampling times $%
(t_{n})_{n\geq0}$. This requires introducing the measure $\mathbb
{P}_{\eta
,(p,n)}\in\mathcal{P}(E_{n}^{\prime})$ defined for any pair of
integers $%
0\leq p\leq n$ and any $\eta\in\mathcal{P}(E_{p}^{\prime})$ by
%
\begin{eqnarray}\label{eq:measureetathenmarkov}
\mathbb{P}_{\eta,(p,n)}(dx_{0:n})& :=&\eta
(dx_{0:p})M_{p+1}(x_{p},dx_{p+1})\cdots M_{n}(x_{n-1},dx_{n}) \nonumber\\
& \in&\mathcal{P} \bigl( E_{p}^{\prime}\times ( E_{p+1}\times\cdots
\times E_{n} ) \bigr) =\mathcal{P}(E_{n}^{\prime}),
\end{eqnarray}
where $dx_{0:n}$ denotes an infinitesimal neighborhood of a path
sequence $%
x_{0:n}\in E_{n}^{\prime}$.

\label{adpttt} Given $\mathcal{H}_{p,q}^{(n)}$, with $n\geq0$ and
$0\leq
p\leq q$, we define an increasing sequence of deterministic time steps $
(t_{n})_{n\geq0}$ and a flow of Feynman--Kac measures $(\eta
_{n},\widehat{%
\eta}_{n})$ by induction as follows. Suppose that the resampling time $
t_{n} $ is defined as well as $(\eta_{n},\widehat{\eta}_{n})\in
\mathcal{P}%
(E_{t_{n}}^{\prime})^{2}$. The resampling time $t_{n+1}$ is defined as the
first time ($s>t_{n}$) the quantity $\mathcal{H}_{t_{n},s}^{(n)} (
\mathbb{%
P}_{\widehat{\eta}_{n},(t_{n},s)} ) $ hits the set $%
I_{n}=[a_{n},\infty)$; that is, $t_{n+1}:=\inf{ \{ t_{n}<s\dvt\mathcal
{H}%
_{t_{n},s}^{(n)} ( \mathbb{P}_{\widehat{\eta}_{n},(t_{n},s)} ) \in
I_{n} \} }$. Given $t_{n+1}$, we set
%
\begin{equation}\label{propprod}
\eta_{n+1}=\mathbb{P}_{\widehat{\eta}_{n},(t_{n},t_{n+1})}
\quad \mbox{and}\quad  \widehat{\eta}_{n+1}=\Psi_{n+1} ( \eta_{n+1} )
\end{equation}
with the Boltzmann--Gibbs transformation $\Psi_{n+1}$ associated with the
potential function $\mathcal{G}_{n+1}=W_{t_{n},t_{n+1}}$.

By definition of the Markov transition $\mathcal{M}_{n+1}$ of the excursion
model $\mathcal{X}_{n}$ defined in Section~\ref{secexc}, it can be checked
that
%
\begin{equation}\label{propprod2}
\eta_{n+1}=\mathbb{P}_{\widehat{\eta}_{n},(t_{n},t_{n+1})}=\widehat
{\eta}_{n}\mathcal{M}_{n+1}.
\end{equation}
This yields the recursion (\eqref{propprod} and \eqref{propprod2}) $%
\Longrightarrow\eta_{n+1}=\Psi_{n} ( \eta_{n} ) \mathcal{M}%
_{n+1}.$ Hence the flow of measures $\eta_{n}$ and $\widehat{\eta}_{n}$
coincide with the Feynman--Kac flow of distributions defined in~(\ref{qf})
with the Markov chain and potential function $(\mathcal
{X}_{n},\mathcal{G}%
_{n})$ on excursion spaces defined in~(\ref{refxa}) and (\ref
{laformm}). The
SMC\ approximation of these distributions is studied in Section~\ref{adpative}.

\subsection{Some applications}\label{applisec}

In this section, we examine the inductive construction of the deterministic
resampling times $ ( t_{n} ) _{n\geq0}$ introduced in Section~\ref%
{adpttt} for the criteria (\ref{eq:limitingvariance}) and (\ref%
{eq:limitingentropy}).

$\bullet$ \textit{Squared coefficient of variation}. In this case, we have
\[
\mathcal{H}_{t_{n},s}^{(n)} \bigl( \mathbb{P}_{\widehat{\eta}%
_{n},(t_{n},s)} \bigr) =\frac{\mathbb{E}_{n,\widehat{\eta}_{n}} (
W_{t_{n},s} ( X_{t_{n}+1:s} ) ^{2} ) }{\mathbb{E}_{n,\widehat{%
\eta}_{n}} ( W_{t_{n},s} ( X_{t_{n}+1:s} ) ) ^{2}}-1.
\]
The mappings $s\mapsto\mathcal{H}_{t_{n},s}^{(n)} ( \mathbb{P}_{%
\widehat{\eta}_{n},(t_{n},s)} ) $ are generally increasing. One
natural way to control these variances is to choose an interval $%
I_{n}:=[a_{n},\infty)$, with $a_{n}>0$, then
\[
t_{n+1}:=\inf \{t_{n}<s\dvt\mathbb{E}_{n,\widehat{\eta}_{n}} (%
W_{t_{n},s} ( X_{t_{n}+1:s} ) ^{2} ) \geq [1+a_{n}%
]{\mathbb{E}_{n,\widehat{\eta}_{n}} (W_{t_{n},s} (
X_{t_{n}+1:s} ) )^{2}} \}.
\]

$\bullet$ \textit{Entropy}. This criterion allows us to control an
entropy-like distance between the free motion trajectories and the weighted
Feynman--Kac measures. To be more precise, set
\[
\mathcal{H}_{t_{n},s}^{(n)} \bigl( \mathbb{P}_{\widehat{\eta}%
_{n},(t_{n},s)} \bigr) =\operatorname{Ent} \bigl( \mathbb{P}_{\widehat{\eta}%
_{n},(t_{n},s)}|\mathbb{Q}_{\widehat{\eta}_{n},(t_{n},s)} \bigr) =-%
\mathbb{E}_{n,\widehat{\eta}_{n}} ( \log W_{t_{n},s} (
X_{t_{n}+1:s} ) )
\]
with the weighted measures $\mathbb{Q}_{\eta,(p,n)}$ defined by
\[
\mathbb{Q}_{\eta,(p,n)}(dx_{0:n})=\mathbb{P}_{\eta
,(p,n)}(dx_{0:n})\times
W_{p,n} ( x_{p+1:n} ) .
\]
If we choose an interval $I_{n}:=[a_{n},\infty)$, with $a_{n}>0$, then the
resampling time $t_{n+1}$ coincides with the first time the entropy
distance goes above the level $a_{n}$; that is,
\[
t_{n+1}:=\inf{ \bigl\{ t_{n}<s\dvt\operatorname{Ent} \bigl( \mathbb{P}_{\widehat{%
\eta}_{n},(t_{n},s)}|\mathbb{Q}_{\widehat{\eta}_{n},(t_{n},s)} \bigr)
\geq a_{n} \bigr\} }.
\]

\section{Convergence analysis of the reference SMC\ algorithm}
\label{adpative}


%

\subsection{A reference SMC\ algorithm}
\label{refsec}

The SMC interpretation of the evolution equation (\ref{predeq}) is the Markov chain
\[
\mathcal{X}_{n}^{(N)}= \bigl( \mathcal{X}_{n}^{(N,1)},\mathcal{X}%
_{n}^{(N,2)},\ldots,\mathcal{X}_{n}^{(N,N)} \bigr) \in S_{n}^{N}
\]
with elementary transitions
%
\begin{equation}\label{elemnta}
\mathbb{P} \bigl( \mathcal{X}_{n+1}^{(N)}\in B_{1}\times\cdots\times
B_{N}|\mathcal{X}_{n}^{(N)} \bigr) =\int_{B_{1}\times\cdots\times
B_{N}}\prod_{i=1}^{N}\mathcal{K}_{n+1,\eta_{n}^{N}}\bigl(\mathcal{X}%
_{n}^{(N,i)},dx_{n+1}^{i}\bigr),
\end{equation}
where $B_{i}\in\mathscr{S}_{n+1}$ for every $i\in\{1,\dots,N\}$ and
%
\begin{equation}\label{empiricalreferenceSMC}
\mathcal{K}_{n+1,\eta_{n}^{N}}=\mathcal{S}_{n,\eta_{n}^{N}}\mathcal
{M}_{n+1} \quad \mbox{and}\quad  \eta_{n}^{N}(\cdot):=\frac{1}{N}%
\sum_{j=1}^{N}\delta_{\mathcal{X}_{n}^{(N,j)}}(\cdot).
\end{equation}
This integral decomposition shows that the SMC\ algorithm has a similar
updating/prediction nature as the one of the `limiting' Feynman--Kac model.
More precisely, the deterministic two-step updating/prediction transitions
in distribution spaces
%
\begin{equation}
\eta_{n}\onrarrow{\mathcal{S}_{n,\eta_{n}}}\widehat{\eta}_{n}=\eta
_{n}\mathcal{S}_{n,\eta_{n}}=\Psi_{n}(\eta_{n})\onrarrow{\mathcal
{M}_{n+1}}%
\eta_{n+1}=\widehat{\eta}_{n}\mathcal{M}_{n+1} \label{flowli}
\end{equation}
have been replaced by a two-step resampling/mutation transition in a product
space
%
\begin{equation}
\mathcal{X}_{n}^{(N)}\in S_{_{n}}^{N}\onrarrow{\mathrm{resampling}}\widehat{\mathcal{X}}_{n}^{(N)}\in
S_{_{n}}^{N}\onrarrow{\mathrm{mutation}}\mathcal{X}_{n+1}^{(N)}\in S_{_{n+1}}^{N}. \label{flowlia}
\end{equation}
In our context, the SMC algorithm keeps track of all the paths of the
sampled particles and the corresponding ancestral lines are denoted by $
\widehat{\mathcal{X}}_{n}^{(N,i)}=\widehat{X}_{0:t_{n}}^{(N,i)}$ and
$%
\mathcal{X}_{n}^{(N,i)}=X_{0:t_{n}}^{(N,i)}\in S_{_{n}}$, where we recall
that $S_{n}=E_{t_{n}}^{\prime}$. By definition of the reference Markov
model~$\mathcal{X}_{n}$ given in (\ref{refxa}), every path particle~$%
\mathcal{X}_{n+1}^{(N,i)}\in S_{_{n+1}}$ keeps track of the selected
excursion~$\widehat{\mathcal{X}}_{n}^{(N,i)}\in S_{_{n}}$ and it evolves
from its terminal state $\widehat{X}_{t_{n},t_{n}}^{(N,i)}$ with $%
(t_{n+1}-t_{n})$ elementary moves using the Markov transition $
M_{t_{n}+1:t_{n+1}}$. More formally, we have that
\[
\mathcal{X}_{n+1}^{(N,i)}= \bigl( \widehat{X}%
_{0:t_{n}}^{(N,i)},X_{t_{n}+1:t_{n+1}}^{(N,i)} \bigr) = \bigl( \widehat{%
\mathcal{X}}_{n}^{(N,i)},X_{t_{n}+1:t_{n+1}}^{(N,i)} \bigr) .
\]
From this discussion, it is worth mentioning a further convention that the
particle empirical measures $\eta_{n+1}^{N}(\cdot)=\frac{1}{N}%
\sum_{i=1}^{N}\delta_{ ( \widehat{\mathcal{X}}%
_{n}^{(N,i)},X_{t_{n}+1:t_{n+1}}^{(N,i)} ) }(\cdot)$ are the terminal
values at time $s=t_{n+1}$ of the flow of random measures
%
\begin{equation}\label{empiref}
t_{n}\leq s\leq t_{n+1}\mapsto\mathbb{P}_{\widehat{\eta}%
_{n}^{N},(t_{n},s)}^{N}(\cdot)=\frac{1}{N}\sum_{i=1}^{N}\delta_{ (
\widehat{\mathcal{X}}_{n}^{(N,i)},X_{t_{n}+1:s}^{(N,i)} ) }(\cdot).
\end{equation}

\subsection{Concentration analysis}\vspace*{-3pt}
\label{conccsecc}

\subsubsection{Introduction}\vspace*{-2pt}
\label{introconcc}

This section is concerned with the concentration analysis
of the empirical measures $\eta_{n}^{N}$ associated with (\ref%
{empiricalreferenceSMC}) around their limiting values $\eta_{n}$
defined in
(\ref{qf}). Our concentration estimates are expressed in terms of\vspace*{-2pt}
\[
q_{p,n}=\sup_{(x,y)\in S_{p}^{2}}\frac{\mathcal
{Q}_{p,n}(1)(x)}{\mathcal{Q}%
_{p,n}(1)(y)} \quad \mbox{and}\quad  \beta(\mathcal{P}_{p,n}):=\sup_{f\in%
\operatorname{Osc}_{1}(S_{n})}\operatorname{osc}(\mathcal{P}_{p,n}(f))\vspace*{-3pt}
\]
with $\mathcal{Q}_{p,n}$ as in (\ref{qpnref}) and $\mathcal
{P}_{p,n}$ in
equation~(\ref{qnormref}). These parameters can be expressed in terms
of the
mixing properties of the Markov transitions $\mathcal{M}_{n}$; see
\cite{fk}, %
Chapter 4. Under appropriate mixing type properties we can prove that
the series $\sum_{p=0}^{n}q_{p,n}^{\alpha}\beta(\mathcal{P}_{p,n})$ is
uniformly bounded w.r.t. the final time horizon $n$ for any parameter $%
\alpha\geq0$. Most of the results presented in this section are expressed
in terms of these series. As a result, these non-asymptotic results can be
converted into time uniform convergence results. To get a flavor of these
uniform estimates, assume that the Markov transitions $\mathcal{M}_{k}$
satisfy the following regularity property.

$(\mathcal{M})_{m}$ There exists an $m\in\mathbb{N}$ and a
sequence $(\delta_{p})_{p\geq0}\in ( 0,1 ) ^{\mathbb{N}}$ such
that\vspace*{-2pt}
\[
\forall p\geq0,\ \forall(x,y)\in S_{p}^{2} \qquad \mathcal{M}%
_{p,p+m}(x,\cdot)\geq\delta_{p}\mathcal{M}_{p,p+m}(y,\cdot)\vspace*{-3pt}
\]
with $\mathcal{M}_{p,p+m}:=\mathcal{M}_{p+1}\mathcal{M}_{p+2}\cdots
\mathcal{M%
}_{p+m}$.

We also introduce the following quantities:\vspace*{-2pt}
\[
\forall k\leq l\qquad  r_{k,l}:=\sup {\prod_{k\leq p<l}\frac{\mathcal
{G}%
_{p}(x_{p})}{\mathcal{G}_{p}(y_{p})}} \qquad \biggl( \leq\prod_{k\leq
p<l}q_{p} \biggr)\vspace*{-3pt}
\]
with the collection of constants $(q_{n})_{\geq1}$ introduced in (\ref
{condG}%
). In the above displayed formula, the supremum is taken over all admissible
pairs of paths with elementary
transitions~$\mathcal{M}_{p}$.

Under the condition $(\mathcal{M})_{m}$ we have for any $n\geq m\geq
1$, and
$p\geq1$,\vspace*{-2pt}
%
\begin{equation} \label{prodd1}
q_{p,p+n}\leq\delta_{p}^{-1}r_{p,p+m}\quad  \mbox{and} \quad \beta(\mathcal{P}_{p,p+n})\leq\prod_{k=0}^{\lfloor n/m\rfloor
-1} \bigl( 1-\delta_{p+km+1}^{2}r_{p+km+1,p+(k+1)m}^{-1} \bigr) .\vspace*{-3pt}
\end{equation}
The proof of these estimates relies on semigroup techniques; see \cite
{fk}, %
Chapter 4, for details. Several contraction inequalities can be deduced
from these results. To understand this more closely, assume that
$(\mathcal{M%
})_{m}$ is satisfied with $m=1$, $\delta=\bigwedge_{n}{\delta_{n}}>0$
and $%
q=\bigvee_{n\geq1}q_{n}$. In this case, $q_{p,p+n}\leq\delta^{-1}q$
and $%
\beta(\mathcal{P}_{p,p+n})\leq ( 1-\delta^{2} ) ^{n}$ imply that\vspace*{-2pt}
\[
\forall\alpha\geq0\qquad  \sum_{p=0}^{n}q_{p,n}^{\alpha}\beta(\mathcal{P
}_{p,n})\leq q^{\alpha}/\delta^{(2+\alpha)}.\vspace*{-3pt}
\]

More generally, assume $(\mathcal{M})_{m}$ is satisfied for some
$m\geq1$
and that the parameters $\delta_{p}$ and $r_{k,l}$ are such that\vspace*{-2pt}
%
\begin{equation}\label{hyprref}
\bigwedge_{p}\delta_{p}:=\delta>0 ,\qquad \bigvee_{p}r_{p,p+m}:=\overline{r}<\infty
\quad \mbox{and}\quad
\bigvee_{p}r_{p+1,p+m}:=\underline{r}<\infty.\vspace*{-3pt}\vadjust{\goodbreak}
\end{equation}
In this situation, $q_{p,p+n}\leq\delta^{-1}\overline{r}$ and $\beta
(\mathcal{P}_{p,p+n})\leq ( 1-\delta^{2}\underline{r}^{-1} ) ^{\lfloor
n/m\rfloor}$ and therefore
%
\begin{equation}
\forall\alpha\geq0\qquad  \sum_{p=0}^{n}q_{p,n}^{\alpha}\beta(\mathcal{P
}_{p,n})\leq m\underline{r}\overline{r}^{\alpha}/\delta^{(2+\alpha)}.
\label{sumcv}
\end{equation}
See \cite{fk}, Chapter 3, for a discussion of when $(\mathcal{M})_{m}$
holds. We also mention that this mixing condition is never met for
$\mathcal{%
X}_{n}=(X_{p})_{0\leq p\leq t_{n}}$ on $S_{n}=E_{t_{n}}^{\prime}$ discussed
in Section~\ref{secexc}. Nevertheless, under appropriate conditions on the
Markov transitions $M_{k}$, it is satisfied for the time marginal model
associated with the excursion valued Markov chain model on $%
\prod_{t_{n-1}<p\leq t_{n}}E_{p}$. 
For instance, if $\forall k\geq1$, $\forall(x,y)\in ( E_{k} )
^{2}$, $M_{k}(x,\cdot)\geq\delta^{\prime}M_{k}(y,%
\cdot)$ for some $\delta^{\prime}>0$, then condition $(\mathcal
{M})_{m}$ is
met with $m=1$ and $\delta_{p}=\delta^{\prime}$.

\subsubsection{Some $\mathbb{L}_{m}$-mean error bounds}
\label{lmbounds}

At this point, it is convenient to observe that the local sampling errors
induced by the mean field particle model are expressed in terms of the
collection of local random field models defined below.

\begin{defi}
For any $n\geq0$ and any $N\geq1$, let $V_{n}^{N}$ be the collection of
random fields defined by the following stochastic perturbation formulae
%
\begin{equation} \label{defVN}
\eta^{N}_{n}=\eta_{n-1}^{N} \mathcal{K }_{n,\eta_{n-1}^{N}}+\frac
{1}{\sqrt{N}%
}V^{N}_{n}\qquad  \bigl( \Longleftrightarrow V^{N}_{n}:=\sqrt{N} [
\eta^{N}_{n}-\eta_{n-1}^{N} \mathcal{K }_{n,\eta_{n-1}^{N}} ] \bigr) .
\end{equation}
For $n=0$, the conventions $\mathcal{K}_{0,\eta_{-1}^{N}}(x,\mathrm{d}y)=\eta_{0}(\mathrm{d}y)$
and $\eta_{-1}^{N} \mathcal{K }_{0,\eta_{-1}^{N}}=\eta_{0} $ are adopted.
\end{defi}

In order to quantify high-order $\mathbb{L}_{m}$-mean errors we need the
following Khinchine type inequality for martingales with symmetric and
independent increments. This is a well-known result.

\begin{lem}[(Khinchine's inequality)]
\label{lemP} Let $L_{n}^{\Delta}:=\sum_{0\leq p\leq n}\Delta_{p}$
be a
real-valued martingale with symmetric and independent increments
$(\Delta
_{n})_{n\geq0}$. For any integer $m\geq1$ and any $n\geq0$, we have%
%
\begin{equation} \label{bdg}
\mathbb{E} ( \vert L_{n}^{\Delta} \vert^{m} ) ^{{1%
}/{m}}\leq b(m)\mathbb{E} ( [ L^{\Delta} ] _{n}^{m^{\prime
}/2} ) ^{{1}/{m^{\prime}}} \qquad \mbox{with } [ L^{\Delta}%
] _{n}:=\sum_{0\leq p\leq n}\Delta_{p}^{2},
\end{equation}
where $m^{\prime}$ stands for the smallest even integer $m^{\prime
}\geq m$
and $(b(m))_{m\geq1}$ is the collection of constants given below:
%
\begin{equation}\label{collec}
b(2m)^{2m}:=(2m)_{m}2^{-m} \quad \mbox{and}\quad  b(2m+1)^{2m+1}:=\frac{%
(2m+1)_{(m+1)}}{\sqrt{m+1/2}}2^{-(m+1/2)}
\end{equation}
with $(2m)_{m}=(2m)!/(2m-m)!$.
\end{lem}

\begin{prop}
\label{lmvp} For any $N\geq1$, $m\geq1$, $n\geq0$ and any test
function $%
f_{n}\in\mathcal{B }_{b}(S_{n})$ we have the almost sure estimate
%
\begin{equation} \label{reflm}
\mathbb{E} \bigl( | V^{N}_{n}(f_{n}) | ^{m} | \mathcal{F}%
^{(N)}_{n-1} \bigr) ^{{1}/{m}}\leq b(m)\operatorname{osc}(f_{n}),
\end{equation}
where $(\mathcal{F}_n^{(N)})_{n\geq0}$ is the filtration generated by
the $%
N $-particle system.
\end{prop}

\begin{pf}
By construction, we have
\begin{eqnarray*}
V_{n}^{N}(f_{n}) &\hspace*{3pt}=&\sum_{i=1}^{N}\Delta_{n,i}^{(N)}(f_{n}), \\
\Delta_{n,i}^{(N)}(f_{n})&:= &\frac{1}{\sqrt{N}}\bigl [ f_{n}\bigl(\mathcal{X}
_{n}^{(N,i)}\bigr)-\mathcal{K}_{n,\eta_{n-1}^{N}}(f_{n})\bigl(\mathcal{X}%
_{n-1}^{(N,i)}\bigr) \bigr] .
\end{eqnarray*}
Given $\mathcal{X}_{n-1}^{(N)}$, let $(\mathcal
{Y}_{n}^{(N,i)})_{1\leq i\leq
N}$ be an independent copy of $(\mathcal{X}_{n}^{(N,i)})_{1\leq i\leq N}$.
It can be checked that
\[
\Delta_{n,i}^{(N)}(f_{n})=\mathbb{E} \biggl( \frac{1}{\sqrt{N}} \bigl[ f_{n}\bigl(%
\mathcal{X}_{n}^{(N,i)}\bigr)-f_{n}\bigl(\mathcal{Y}_{n}^{(N,i)}\bigr) \bigr] \big|\mathcal
{F}%
_{n}^{(N)} \biggr).
\]
This yields the formula $V_{n}^{N}(f_{n})=\mathbb{E}(L_{n,N}^{(N)}(f_{n})|
\mathcal{F}_{n}^{(N)}),$ where $L_{n,N}^{(N)}(f_{n})$ is the terminal value
of the martingale sequence defined by
\[
i\in\{1,\ldots,N\}\mapsto L_{n,i}^{(N)}(f_{n}):=\frac{1}{\sqrt{N}}%
\sum_{j=1}^{i} \bigl[ f_{n}\bigl(\mathcal
{X}_{n}^{(N,j)}\bigr)-f_{n}\bigl(\mathcal{Y%
}_{n}^{(N,j)}\bigr) \bigr] .
\]
Then as
\begin{eqnarray*}
\mathbb{E}\bigl(|V_{n}^{N}(f_{n})|^{m}|\mathcal{F}_{n-1}^{(N)}\bigr)^{1/m}
&=&\mathbb{E%
} \bigl(\bigl|\mathbb{E} \bigl(L_{n,N}^{(N)}(f_{n})|\mathcal{F}_{n}^{(N)} \bigr)\bigr|^{m}\big|%
\mathcal{F}_{n-1}^{(N)} \bigr)^{1/m} \\
&\leq&\mathbb{E} \bigl(\bigl|L_{n,N}^{(N)}(f_{n})\bigr|^{m}\big|\mathcal{F}_{n-1}^{(N)}%
\bigr)^{1/m},
\end{eqnarray*}
one may apply Khinchine's inequality to conclude.
\end{pf}

The proof of the following lemma is rather technical and is provided in the
\hyperref[app]{Appendix}.

\begin{lem}
\label{keydec} For any $0\leq p\leq n$, any $\eta,\mu\in\mathcal{P}
(S_{p}) $ and any $f_{n}\in\operatorname{Osc}_{1}(S_{n})$, we have
the first-order decomposition for the nonlinear semigroup $\Phi_{p,n}$
defined in (%
\ref{phiQ}):
\[
[ \Phi_{p,n}(\mu)-\Phi_{p,n}(\eta) ] (f_{n})=2q_{p,n}\beta(%
\mathcal{P}_{p,n}) [ \mu-\eta ] ( \mathcal{U}_{p,n,\eta
}(f_{n}) ) +\mathcal{R}_{p,n}(\mu,\eta)(f_{n}),
\]
where
\[
\vert\mathcal{R}_{p,n}(\mu,\eta)(f_{n}) \vert\leq
4q_{p,n}^{3}\beta(\mathcal{P}_{p,n}) \vert [ \mu-\eta ]
( \mathcal{V}_{p,n,\eta}(f) ) \vert\times \vert [
\mu-\eta ] ( \mathcal{W}_{p,n,\eta}(f_{n}) ) \vert
\]
with $\mathcal{U}_{p,n,\eta}(f),\mathcal{V}_{p,n,\eta}(f),\mathcal
{W}%
_{p,n,\eta}(f)$ a collection of functions in $\operatorname{Osc}_{1}(S_{p})$
whose values only depend on the parameters $(p,n,\eta)$.
\end{lem}

We now present a bias estimate and some $\mathbb{L}_{m}$ bounds of
independent interest.

\begin{theo}
\label{interm} For any $n\geq0$, $f_{n}\in\operatorname
{Osc}_{1}(S_{n})$ and any
$N\geq1$,
\[
N \vert\mathbb{E} ( \eta_{n}^{N}(f_{n}) )
-\eta_{n}(f_{n}) \vert\leq\sigma_{1,n} \qquad \mbox{with } \sigma
_{1,n}:=4\sum_{p=0}^{n}q_{p,n}^{3}\beta(\mathcal{P}_{p,n}).
\]
In addition, for any $m\geq1$ we have
\[
\sqrt{N}\mathbb{E} \bigl( \vert [ \eta_{n}^{N}-\eta_{n} ]
(f_{n}) \vert^{m} \bigr) ^{{1}/{m}}\leq\frac{1}{\sqrt{N}}%
b(2m)^{2}\sigma_{1,n}+b(m)\sigma_{2,n}
\]
with $\sigma_{2,n}:=2\sum_{p=0}^{n}q_{p,n}\beta(\mathcal{P}_{p,n}). $
\end{theo}

\begin{pf}
Using Lemma \ref{keydec}, we have the telescoping sum decomposition
\begin{eqnarray*}
W_{n}^{N}& :=&\sqrt{N} [ \eta_{n}^{N}-\eta_{n} ] \\
& \hspace*{3pt}=&\sqrt{N}\sum_{p=0}^{n} [ \Phi_{p,n}(\eta_{p}^{N})-\Phi
_{p,n} ( \Phi_{p}(\eta_{p-1}^{N}) ) ] =\mathcal{I}_{n}^{N}+%
\mathcal{J}_{n}^{N}
\end{eqnarray*}
with $\eta_{-1}^{N}(f):=f$ and the pair of random measures $(\mathcal
{I}%
_{n}^{N},\mathcal{J}_{n}^{N})$ given for any $f_{n}\in\operatorname
{Osc}%
_{1}(S_{n})$ by
\begin{eqnarray*}
\mathcal{I}_{n}^{N}(f_{n})& :=&2\sum_{p=0}^{n}q_{p,n}\beta(\mathcal
{P}%
_{p,n})V_{p}^{N} \bigl( \mathcal{U}_{p,n,\Phi_{p}(\eta
_{p-1}^{N})}(f_{n}) \bigr), \\
\mathcal{J}_{n}^{N}(f_{n})& :=&\sqrt{N}\sum_{p=0}^{n}\mathcal
{R}_{p,n} (
\eta_{p}^{N},\Phi_{p}(\eta_{p-1}^{N}) ) (f_{n}).
\end{eqnarray*}
Now, observe that
%
\begin{equation}
\mathbb{E} ( W_{n}^{N}(f_{n}) ) =\mathbb{E}(\mathcal{J}%
_{n}^{N}(f_{n})). \label{eq:wneqin}
\end{equation}
Using Proposition~\ref{lmvp}, for any $f_{n}\in\operatorname{Osc}_{1}(S_{n})$
it can be checked that $\mathbb{E} ( \vert\mathcal{I}%
_{n}^{N}(f_{n}) \vert^{m} ) ^{{1}/{m}}\leq b(m)\sigma
_{2,n}. $ In a similar way, we find that
%
\begin{equation}\label{LmI}
\sqrt{N}\mathbb{E} ( \vert\mathcal{J}_{n}^{N}(f_{n}) \vert
^{m} ) ^{{1}/{m}}\leq b(2m)^{2}\sigma_{1,n}.
\end{equation}
The first part of the proof then follows from \eqref{eq:wneqin} and %
\eqref{LmI}; the remainder of the proof is now clear.
\end{pf}

\subsubsection{A concentration theorem}
\label{conctheo}

The following concentration theorem is the main result of this section.
\begin{theo}
\label{theoexpoc} For any $n\geq0$, $f_{n}\in\operatorname
{Osc}_{1}(S_{n})$, $%
N\geq1$ and any $0\leq\epsilon\leq1/2$,
%
\begin{equation}\label{conccineg}
\mathbb{P} \bigl( \vert [ \eta_{n}^{N}-\eta_{n} ]
(f_{n}) \vert\geq\epsilon\bigr) \leq6\exp{ \biggl( -\frac{%
N\epsilon^{2}}{8\sigma_{1,n}} \biggr), }
\end{equation}
where the constant $\sigma_{1,n}$ is as in Theorem~\ref{interm}.

In addition, suppose $(\mathcal{M})_{m}$ is satisfied for some $m\geq
1$ and
condition (\ref{hyprref}) holds true for some $\delta>0$ and some finite
constants $(\underline{r},\overline{r})$. In this situation, for any value
of the time parameter $n$, for any $f_{n}\in\operatorname{Osc}_{1}(S_{n})$,
$N\geq1$ and for any $\rho\in ( 0,1 ), $ the probability that
\[
\vert [ \eta_{n}^{N}-\eta_{n} ] (f_{n}) \vert\leq
\frac{4\overline{r}}{\delta^{2}}\sqrt{\frac{2m\underline
{r}\overline{r}}{%
N\delta}\log{ \biggl( \frac{6}{\rho} \biggr) }}
\]
is greater than $(1-\rho)$.
\end{theo}

\begin{pf}
We use the same notation as in the proof of Theorem~\ref{interm}. Recall
that $b(2m)^{2m}=\mathbb{E}(X^{2m})$ for every centered Gaussian random
variable with $\mathbb{E}(X^{2})=1$ and
\[
\forall s\in\lbrack0,1/2)\qquad  \mathbb{E}(\exp{ \{ sX^{2} \} }%
)=\sum_{m\geq0}\frac{s^{m}}{m!}b(2m)^{2m}=\frac{1}{\sqrt{1-2s}}.
\]
Using (\ref{LmI}), for any $f_{n}\in\operatorname{Osc}_{1}(S_{n})$
and $0\leq
s<1/(2\sigma_{1,n}),$ it follows that
%
\begin{equation}\label{eqs1n}
\mathbb{E} \bigl( \exp{ \bigl\{ s\sqrt{N}\mathcal{J}_{n}^{N}(f_{n}) \bigr\} }%
\bigr) \leq\sum_{m\geq0}\frac{(s\sigma_{1,n})^{m}}{m!}b(2m)^{2m}=\frac
{1%
}{\sqrt{1-2s\sigma_{1,n}}}.
\end{equation}
To simplify the presentation, set
\[
f_{p,n}^{N}:=\mathcal{U}_{p,n,\Phi_{p}(\eta_{p-1}^{N})}(f_{n}) %
\quad \mbox{and}\quad  \alpha_{p,n}:=2q_{p,n}\beta(\mathcal{P}_{p,n}),
\]
where $\mathcal{U}_{p,n,\eta} ( \cdot ) $ was introduced in Lemma %
\ref{keydec}. By the definition of $V_{p}^{N}$
\[
V_{p}^{N} ( f_{p,n}^{N} ) =\frac{1}{\sqrt{N}}\sum_{i=1}^{N} \bigl(
f_{p,n}^{N}\bigl(\mathcal{X}_{p}^{(N,i)}\bigr)-\mathcal{K}_{p,\eta
_{p-1}^{N}}(f_{p,n}^{N})\bigl(\mathcal{X}_{p-1}^{(N,i)}\bigr) \bigr) .
\]
Recalling that $\mathbb{E}(\mathrm{e}^{tX})\leq \mathrm{e}^{t^{2}c^{2}/2}$ for every
real-valued and centered random variable $X$ with \mbox{$|X|\leq c$} (e.g.,
\cite{fk}, %
Lemma 7.3.1), we prove that
\begin{eqnarray*}
&&\mathbb{E} \bigl( \exp{ \{ t\alpha_{p,n}V_{p}^{N} (
f_{p,n}^{N} ) \} } \vert\mathcal{X}_{p-1}^{(N)}
\bigr) \\[3pt]
&&\quad =\prod_{i=1}^{N} \int_{S_{p}}\mathcal{K}_{p,\eta
_{p-1}^{N}} \bigl( \mathcal{X}_{p-1}^{(N,i)},\mathrm{d}x \bigr) \mathrm{e}^{(t\alpha
_{p,n}/{\sqrt{N}}) ( f_{p,n}^{N}(x)-\mathcal{K}_{p,\eta
_{p-1}^{N}}(f_{p,n}^{N})(\mathcal{X}_{p-1}^{(N,i)}) ) }\leq\exp{ ({%
t^{2}\alpha_{p,n}^{2}}/{2} )}.%
\end{eqnarray*}
Iterating the argument, we find that
%
\begin{equation}\label{eqs22}
\mathbb{E} \bigl(\mathrm{e}^{t\mathcal{I}_{n}^{N}(f_{n})} \bigr)=\mathbb{E} \Biggl(\exp{%
\Biggl\{t\sum_{p=0}^{n}\alpha_{p,n}V_{p}^{N} (f_{p,n}^{N} ) \Biggr\}} \Biggr)%
\leq\exp{ \biggl(\frac{t^{2}\sigma_{n}^{2}}{2} \biggr)}
\end{equation}
with $\sigma_{n}^{2}:=4\sum_{p=0}^{n}q_{p,n}^{2}\beta(\mathcal{P}%
_{p,n})^{2}.$

From these upper bounds, the proof of the exponential estimates now follows
standard arguments. Indeed, for any $0\leq s<{1}/{(2\sigma_{1,n})}$
and any
$\epsilon>0$, by (\ref{eqs1n}) we have
\[
\mathbb{P} \bigl( \sqrt{N}\mathcal{J}_{n}^{N}(f_{n})\geq\epsilon \bigr)
\leq\frac{1}{\sqrt{1-2s\sigma_{1,n}}}\exp{ \{ -\epsilon s \} }.
\]
Replacing $\epsilon$ by $\epsilon N$ and choosing $s={3}/{(8\sigma_{1,n})}$
yields
\[
\mathbb{P} \bigl( \mathcal{J}_{n}^{N}(f_{n})/\sqrt{N}\geq\epsilon \bigr)
\leq2\exp{ \{ -\epsilon N/(3\sigma_{1,n}) \} }.
\]
To estimate the probability tails of $\mathcal{I}_{n}^{N}(f_{n})$, we
use (%
\ref{eqs22}) and the fact that $\epsilon>0$ and $t\geq0$
\[
\mathbb{P} \bigl( \mathcal{I}_{n}^{N}(f_{n})\geq\epsilon \bigr) \leq\exp{%
\biggl\{ - \biggl( \epsilon t-\frac{t^{2}}{2}\sigma_{n}^{2} \biggr) \biggr\} }.
\]
Now, choosing $t=\epsilon/\sigma_{n}^{2}$ and replacing $\epsilon$
by $%
\sqrt{N}\epsilon$, we obtain
\[
\forall\epsilon>0\qquad  \mathbb{P} \bigl( \mathcal{I}_{n}^{N}(f)/\sqrt{N}%
\geq\epsilon \bigr) \leq\exp{ \biggl( -\frac{N\epsilon^{2}}{2\sigma
_{n}^{2}} \biggr) }.
\]
Using the decomposition
\[
[ \eta_{n}^{N}-\eta_{n} ] =\mathcal{I}_{n}^{N}/\sqrt{N}+\mathcal{%
J}_{n}^{N}/\sqrt{N}
\]
we find that for any parameter $\alpha\in\lbrack0,1]$
\[
\mathbb{P} \bigl( [ \eta_{n}^{N}-\eta_{n} ] (f_{n})\geq\epsilon
\bigr) \leq\mathbb{P} \bigl( \mathcal{I}_{n}^{N}(f_{n})/\sqrt{N}\geq
\alpha\epsilon \bigr) +\mathbb{P} \bigl( \mathcal{J}_{n}^{N}(f_{n})/\sqrt{N%
}\geq(1-\alpha)\epsilon \bigr) .
\]
From previous calculations,
%
\begin{equation}\label{optim}
\mathbb{P} \bigl( [ \eta_{n}^{N}-\eta_{n} ] (f_{n})\geq\epsilon
\bigr) \leq\exp{ \biggl( -\frac{N\epsilon^{2}\alpha^{2}}{2\sigma_{n}^{2}}%
\biggr) }+2\exp{ \biggl( -\frac{N\epsilon(1-\alpha)}{3\sigma_{1,n}}%
\biggr) }.
\end{equation}
Now, choose $\alpha=(1-\epsilon)(\geq1/2)$, then $\alpha^{2}\geq
1/4$ and
\begin{eqnarray*}
\mathbb{P} \bigl( [ \eta_{n}^{N}-\eta_{n} ] (f_{n})\geq\epsilon
\bigr) &\leq&\exp{ \biggl( -\frac{N\epsilon^{2}}{8\sigma_{n}^{2}} \biggr) }%
+2\exp{ \biggl( -\frac{N\epsilon^{2}}{3\sigma_{1,n}} \biggr) } \\
&\leq&3\exp{ \biggl( -\frac{N\epsilon^{2}}{8 ( \sigma_{1,n}\vee
\sigma_{n}^{2} ) } \biggr) }.
\end{eqnarray*}
It remains to observe that $q_{p,n}\geq1$ and $\beta(\mathcal{P}%
_{p,n})\leq1\Longrightarrow\sigma_{n}^{2}\leq\sigma_{1,n}$ and
\begin{eqnarray*}
\vert [ \eta_{n}^{N}-\eta_{n} ] (f_{n}) \vert\geq
\epsilon&\quad  \Longleftrightarrow\quad  &[ \eta_{n}^{N}-\eta_{n} ]
(f_{n})\geq\epsilon \quad \mbox{or}\quad  [ \eta_{n}^{N}-\eta_{n}%
] (f_{n})\leq-\epsilon\\
&\quad  \Longleftrightarrow\quad &
 [ \eta_{n}^{N}-\eta_{n} ] (f_{n})\geq
\epsilon \quad \mbox{or}\quad  [ \eta_{n}^{N}-\eta_{n} ]
(-f_{n})\geq\epsilon
\end{eqnarray*}
so that
\[
\mathbb{P} \bigl( \vert [ \eta_{n}^{N}-\eta_{n} ]
(f_{n}) \vert\geq\epsilon \bigr) \leq\mathbb{P} \bigl( [ \eta
_{n}^{N}-\eta_{n} ] (f_{n})\geq\epsilon \bigr) +\mathbb{P} \bigl( %
[ \eta_{n}^{N}-\eta_{n} ] (-f_{n})\geq\epsilon \bigr) .
\]
The end of the proof of (\ref{conccineg}) is now easily completed. We now
assume that the mixing condition $(\mathcal{M})_{m}$ is satisfied for
some $m\geq1$ and condition (\ref{hyprref}) holds true for some $\delta>0$ and
some finite constants $(\underline{r},\overline{r})$. By (\ref
{sumcv}) the
following uniform concentration estimate holds
\[
\sup_{n\geq0} \mathbb{P} \bigl( \vert [ \eta_{n}^{N}-\eta_{n}%
] (f_{n}) \vert\geq\epsilon \bigr) \leq6\exp{ \biggl( -%
\frac{N\epsilon^{2}\delta^{5}}{32m\underline{r}\overline{r}^{3}}%
\biggr) }.
\]
The proof of the theorem is concluded by choosing $\epsilon:=\frac{1}{
\sqrt{N}}\frac{4\overline{r}}{\delta^{2}}\sqrt{\frac{2m\underline
{r}%
\overline{r}}{\delta}\log{\frac{6}{\rho}}}$.
\end{pf}

\begin{rem}
Returning to the end of the proof of Theorem~\ref{theoexpoc}, the
exponential concentration estimates can be marginally improved by
choosing, in (\ref{optim}), the parameter $\alpha=\alpha
_{n}(\epsilon)\in[0,1]$ such that $a_{n}(\epsilon)\alpha
^{2}=b_{n}(1-\alpha)
$, with $a_{n}(\epsilon):=\frac{\epsilon}{2\sigma^{2}_{n}}$, $b_{n}%
=\frac{1}{3\sigma_{1,n}}$ and $\sigma_{n}^{2}:=4\sum
_{p=0}^{n}q_{p,n}^{2}\beta(\mathcal{P}%
_{p,n})^{2}.$ Elementary manipulations yield
\begin{eqnarray*}
\alpha_{n}(\epsilon)&=&\frac{b_{n}}{2a_{n}(\epsilon)} \biggl( \sqrt
{1+\frac{4a_{n}(\epsilon)}{b_{n}}}-1 \biggr) \\
&=&\frac{\sigma_{n}^{2}}%
{3\sigma_{1,n}}\frac{1}{\epsilon} \biggl( \sqrt{1+\frac{6\sigma_{1,n}}%
{\sigma^{2}_{n}}\epsilon}-1 \biggr) (\longrightarrow_{\epsilon\downarrow0}1)
\end{eqnarray*}
and therefore
\[
\forall\epsilon\geq0 \qquad \mathbb{P} \bigl( | (\eta^{N}_{n}-\eta
_{n}) (f) |\geq\epsilon\bigr) \leq6 \exp{ \biggl( -N\frac
{\epsilon^{2}}{2\sigma^{2}_{n}}\alpha^{2}_{n}(\epsilon) \biggr) }.
\]
For small values of $\epsilon$, this bound improves that in Section
7.4.3 of~\cite{fk}, which is of the form
\[
\forall\epsilon\geq0 \qquad \mathbb{P} \bigl(
| (\eta^{N}_{n}-\eta
_{n}) (f) |\geq\epsilon\bigr) \leq \bigl(1+\epsilon\sqrt{N}\bigr)
\exp{ \biggl( -N\frac{\epsilon^{2}}{2\tilde{\sigma}^{2}_{n}} \biggr)}
\]
with
\[
\tilde{\sigma}^{2}_{n}:=4 \Biggl(\sum_{p=0}^{n}q_{p,n}\beta(\mathcal{P}%
_{p,n}) \Biggr)^2\geq\sigma_{n}^{2}.
\]
\end{rem}

\subsection{Approximating the criteria}
\label{concsec}

By construction, the particle occupation measures~$\mathbb{P}%
_{\widehat{\eta}_{n},(t_{n},s)}^{N}$ approximate the measures~$\mathbb{P}_{%
\widehat{\eta}_{n},(t_{n},s)}$ introduced in (\ref
{eq:measureetathenmarkov}%
); that is, in some sense, $\mathbb{P}_{\widehat{\eta}%
_{n}^{N},(t_{n},s)}^{N}$ $\simeq_{N\uparrow\infty}\mathbb
{P}_{\widehat{%
\eta}_{n},(t_{n},s)}$. Conversely, observe that $\mathbb{P}_{\widehat
{\eta}%
_{n}^{N},(t_{n},s)}^{N}$, respectively $\mathbb{P}_{\widehat{\eta}%
_{n},(t_{n},s)}$, are the marginals of the measures~$\eta_{n+1}^{N}$,
respectively $\eta_{n+1}$, w.r.t. the $(s-t_{n})+1$ first coordinates. In
other words, the measures $\mathbb{P}_{\widehat{\eta}%
_{n}^{N},(t_{n},s)}^{N} $, respectively $\mathbb{P}_{\widehat{\eta}%
_{n},(t_{n},s)}$, are the projections of the measures $\eta_{n+1}^{N}$,
respectively $\eta_{n+1}$, on the state space $E_{s}^{\prime
}=E_{t_{n}}^{\prime}\times ( E_{t_{n}+1}\times\cdots\times
E_{s} ) .$

For instance, the following proposition is essentially a direct consequence
of Theorem~\ref{theoexpoc}.

\begin{prop}\label{propreff}
For any $N\geq1$, $n\geq0$, $t_{n}\leq s\leq t_{n+1}$ and
any $\epsilon>0$, the concentration inequality:
\[
\mathbb{P} \bigl( \bigl\vert\mathcal{H}_{t_{n},s}^{(n)} \bigl( \mathbb{P}_{%
\widehat{\eta}_{n}^{N},(t_{n},s)}^{N} \bigr) -\mathcal{H}_{t_{n},s}^{(n)}%
\bigl( \mathbb{P}_{\widehat{\eta}_{n},(t_{n},s)} \bigr) \bigr\vert
\geq\epsilon\bigr) \leq\bigl(1+\epsilon\sqrt{N/2}\bigr)\exp{ \biggl( -\frac{%
N\epsilon^{2}}{c(n)} \biggr) }
\]
holds for some finite constant $c(n)<\infty$ whose values only depend
on the
time parameter. In addition, when the measures $H_{t_{n},s}^{(n)}$ have a
finite support, the concentration inequality
\[
\mathbb{P} \bigl( \bigl\vert\mathcal{H}_{t_{n},s}^{(n)} \bigl( \mathbb{P}_{%
\widehat{\eta}_{n}^{N},(t_{n},s)}^{N} \bigr) -\mathcal{H}_{t_{n},s}^{(n)}%
\bigl( \mathbb{P}_{\widehat{\eta}_{n},(t_{n},s)} \bigr) \bigr\vert
\geq\epsilon\bigr) \leq c_{1}(n)\exp{ \biggl( -\frac{N\epsilon^{2}}{c_{2}(n)}%
\biggr) }
\]
also holds, with a pair of finite constants $c_{1}(n),c_{2}(n)<\infty$.
\end{prop}

\begin{pf}
By \cite{fk}, Theorem 7.4.4, for any $N\geq1$, $p\geq1$, $n\geq0$
and any
test function $f_{n}\in\operatorname{Osc}_{1}(E_{t_{n}}^{\prime})$
\[
\sup_{N\geq1}\sqrt{N}\mathbb{E} \bigl( \vert\eta_{n}^{N}(f_{n})-\eta
_{n}(f_{n}) \vert^{p} \bigr) ^{{1}/{p}}\leq b(p)c(n)
\]
with some finite constant $c(n)<\infty$ and with the collection of
constants $b(p)$ defined in~(\ref{collec}). These estimates clearly imply
that for any $t_{n}\leq s\leq t_{n+1}$, and any test function $h_{n}\in
\operatorname{Osc}_{1}(E_{s}^{\prime}),$
\[
\sup_{N\geq1}\sqrt{N}\mathbb{E} \bigl( \bigl\vert\mathbb{P}_{\widehat{\eta
}_{n}^{N},(t_{n},s)}^{N}(h_{n})-\mathbb{P}_{\widehat{\eta}%
_{n},(t_{n},s)}(h_{n}) \bigr\vert^{p} \bigr) ^{{1}/{p}}\leq b(p)c(n).
\]
Under (\ref{lipsH}) on the criteria type functionals $\mathcal{H}%
_{t_{n},s}^{(n)}$ and using the generalized integral Minkowski inequality,
it can be concluded that
\[
\sup_{N\geq1}\sqrt{N}\mathbb{E} \bigl( \bigl\vert\mathcal{H}%
_{t_{n},s}^{(n)} \bigl( \mathbb{P}_{\widehat{\eta}_{n}^{N},(t_{n},s)}^{N}%
\bigr) -\mathcal{H}_{t_{n},s}^{(n)} \bigl( \mathbb{P}_{\widehat{\eta}%
_{n},(t_{n},s)} \bigr) \bigr\vert^{p} \bigr) ^{{1}/{p}}\leq
b(p)c(n)\delta\bigl(H_{t_{n},s}^{(n)}\bigr).
\]
The proof of the exponential estimate follows exactly the same lines of
arguments as the ones used in the proof of Corollary 7.4.3 in~\cite{fk};
thus it is omitted. The last assertion is a~direct consequence of
Theorem~%
\ref{theoexpoc}.
\end{pf}

\subsection{An online adaptive SMC\ algorithm}

\label{secYa}

The above proposition shows that the functional criteria $\mathcal{H}%
_{t_{n},s}^{(n)}( \mathbb{P}_{\widehat{\eta}_{n},(t_{n},s)}) $ can be
approximated by $\mathcal{H}_{t_{n},s}^{(n)}( \mathbb{P}_{\widehat
{\eta}%
_{n}^{N},(t_{n},s)}^{N}) $, up to an exponentially small probability.
Therefore, as we cannot compute the deterministic resampling times
$(t_{n})$%
, it is necessary to approximate the reference particle model:

\begin{defi}
\label{defonline} The particle systems $\mathcal{Y}^{(N)}=(\mathcal
{Y}%
^{(N,i)})$, $\widehat{\mathcal{Y}}^{(N)}=(\widehat{\mathcal
{Y}}^{(N,i)})$, $%
Y_{s,t}^{(N,i)}$ and~$\widehat{Y}_{s,t}^{(N,i)}$ are defined as
$\mathcal{X}%
^{(N)}=(\mathcal{X}^{(N,i)})$, $\widehat{\mathcal
{X}}^{(N)}=(\widehat{%
\mathcal{X}}^{(N,i)})$, and $X_{s,t}^{(N,i)}$ and $\widehat{X}_{s,t}^{(N,i)}$
by replacing in the inductive construction of the deterministic
sequence $%
(t_{n})_{n\geq0}$ the measures $\mathbb{P}_{\widehat{\eta}_{n},(t_{n},s)}$
by their current $N$-particle approximation measures $\overline
{\mathbb{P}}_{%
\widehat{\eta}_{n}^{N},(t_{n}^{N},s)}^{N}(\cdot):=\frac{1}{N}%
\sum_{i=1}^{N}\delta_{ (\widehat{\mathcal{Y}}%
_{n}^{(N,i)},Y_{t_{n}^{N}+1:s}^{(N,i)} )}(\cdot).$ Here $\widehat
{\eta}%
_{n}^{N}(\cdot)=\frac{1}{N}\sum_{i=1}^{N}\delta_{\widehat{\mathcal
{Y}}%
_{n}^{(N,i)}}(\cdot)$ denotes the updated occupation measure of the
particle system $\widehat{\mathcal{Y}}_{n}^{(N)}$. We also assume
that both
models are constructed in such a way that they coincide on every time
interval $0\leq n\leq m$, once the random times $t_{n}^{N}=t_{n}$, for every
$0\leq n\leq m$.
\end{defi}

It is emphasized that the measures $\overline{\mathbb{P}}_{\widehat
{\eta}%
_{n}^{N},(t_{n}^{N},s)}^{N}$ differ from the reference empirical
measures $%
\mathbb{P}_{\widehat{\eta}_{n}^{N},(t_{n},s)}^{N}$ in (\ref{empiref}).
Indeed, the reference measures $\mathbb{P}_{\widehat{\eta}%
_{n}^{N},(t_{n},s)}^{N}$ are built using the deterministic times $t_{n}$
based on the functional criteria $\mathcal{H}_{t_{n-1},s}^{(n-1)} (
\mathbb{P}_{\widehat{\eta}_{n-1},(t_{n-1},s)} ) $, whilst the empirical
measures $\overline{\mathbb{P}}_{\widehat{\eta}_{n}^{N},(t_{n}^{N},s)}^{N}$
are inductively constructed using random times $t_{n}^{N}$ based on $%
\mathcal{H}_{t_{n-1}^{N},s}^{(n-1)} (
\overline{\mathbb{P}}_{\widehat{\eta
}_{n-1}^{N},(t_{n-1}^{N},s)}^{N} ) $.

By construction, for the pair of functional criteria discussed in
Section~%
\ref{applisec}, we have that $\mathcal{H}_{t_{n}^{N},s}^{(n)} (
\overline{\mathbb{P}}_{\widehat{\eta}_{n}^{N},(t_{n}^{N},s)}^{N} )
=C_{t_{n}^{N},s}^{N}$, where $C_{t_{n}^{N},s}^{N}$ are the empirical
criteria discussed in Section~\ref{sec:adaptivecriteria}.

\section{Asymptotic analysis}
\label{asymptsec}

\subsection{A key approximation lemma}
\label{secevent}

To go one step further in our discussion, it is convenient
to introduce the following collection of events.

\begin{defi}
For any $\delta\in( 0,1 ) $, $m\geq0$, $a_{n}\in\mathbb{R}$ and $%
N\geq1$, we denote by~$\Omega_{m}^{N} ( \delta,\allowbreak(a_{n})_{0\leq n\leq
m} ) $, the collection of events defined by:
\begin{eqnarray*}
\Omega_{m}^{N} ( \delta,(a_{n})_{0\leq n\leq m} ) & := & \bigl\{
\forall0\leq n\leq m,\forall t_{n}\leq s\leq t_{n+1} \\
& &\hphantom{\bigl\{} \bigl\vert\mathcal{H}_{t_{n},s}^{(n)} \bigl( \mathbb{P}_{\widehat{%
\eta}_{n}^{N},(t_{n},s)}^{N} \bigr) -\mathcal{H}_{t_{n},s}^{(n)} \bigl(
\mathbb{P}_{\widehat{\eta}_{n},(t_{n},s)} \bigr) \bigr\vert
\leq\delta\bigl\vert\mathcal{H}_{t_{n},s}^{(n)} \bigl( \mathbb{P}_{%
\widehat{\eta}_{n},(t_{n},s)} \bigr) -a_{n} \bigr\vert \bigr\} .
\end{eqnarray*}
\end{defi}

The proof of the following result is straightforward and hence omitted.

\begin{lem}
\label{lemke} On the event $\Omega_{m}^{N}(\delta,(a_{n})_{0\leq
n\leq m})$,
for any $n\leq m$ and for any $t_{n}\leq s \leq t_{n+1}$, we have 
%
\[
\mathcal{H }^{(n)}_{t_{n},s} \bigl( \mathbb{P}_{\widehat{\eta}%
_{n},(t_{n},s)} \bigr) >a_{n}\quad \Longrightarrow\quad \mathcal{H }^{(n)}_{t_{n},s}\bigl(
\mathbb{P}^{N}_{\widehat{\eta}_{n}^{N},(t_{n},s)}\bigr) >a_{n}.\vspace*{-3pt}
\]
\end{lem}

\begin{prop}
\label{propke} Assume that the threshold parameters $a_{n}$ are chosen so
that $\mathcal{H}_{t_{n},s}^{(n)} ( \mathbb{P}_{\widehat{\eta}%
_{n},(t_{n},s)} ) \not=a_{n}$, for any $n\geq0$. In this situation,
for any $\delta\in( 0,1 ) $, $m\geq0$ and $N\geq1$, we have
\[
\bigcap_{0\leq n\leq m}\{t_{n}^{N}=t_{n}\}\supset\Omega_{m}^{N} (
\delta,(a_{n})_{0\leq n\leq m} ).\vspace*{-3pt}
\]
\end{prop}

\begin{pf}
This result is proved by induction on $m\geq0$. Under our assumptions,
for $%
m=0$ we have $t_{0}^{N}=t_{0}=0$. Thus, by our coupling construction the
pair of particle models coincide up to the time $ ( t_{1}^{N}\wedge
t_{1} ) $. Therefore, we have
\[
\forall s< ( t_{1}^{N}\wedge t_{1} )\qquad  \mathbb{P}_{\widehat{%
\eta}_{0}^{N},(t_{0},s)}^{N}=\overline{\mathbb{P}}_{\widehat{\eta}%
_{0}^{N},(t_{0}^{N},s)}^{N}.\vspace*{-3pt}
\]
By Lemma~\ref{lemke}, on the event $\Omega_{m}^{N} ( \delta
,(a_{n})_{0\leq n\leq m} ) $ we have $t_{1}^{N}=t_{1}$. This proves the
inclusion for $m=0$ and $m=1$. Suppose the result is true at rank $m$. Thus,
on the event $\Omega_{m}^{N} ( \delta,(a_{n})_{0\leq n\leq m} ) $
it is the case that $t_{n}^{N}=t_{n}$, for any $0\leq n\leq m$. By our
coupling construction, the pair of particle models coincide up to $ (
t_{m+1}^{N}\wedge t_{m+1} ) $; that is,
\[
t_{m}^{N}=t_{m} \quad \mbox{and}\quad  \forall s< (
t_{m+1}^{N}\wedge t_{m+1} )\qquad  \mathbb{P}_{\widehat{\eta}%
_{m}^{N},(t_{m},s)}^{N}=\overline{\mathbb{P}}_{\widehat{\eta}%
_{m}^{N},(t_{m}^{N},s)}^{N}.\vspace*{-3pt}
\]
Once again, by Lemma~\ref{lemke}, on the event $\Omega_{m+1}^{N} (
\delta,(a_{n})_{0\leq n\leq m+1} ) $ it also follows that $%
t_{m+1}^{N}=t_{m+1}$.\vspace*{-2pt}
\end{pf}

\subsection{Randomized criteria}\vspace*{-2pt}
\label{randomizedr}

The situation where the threshold parameters coincide with the adaptive
criteria values $\mathcal{H}_{t_{n},s}^{(n)} ( \mathbb{P}_{\widehat
{\eta
}_{n},(t_{n},s)} ) =a_{n}$ cannot be dealt with using our analysis.
This situation is more involved since it requires us to control both the
empirical approximating criteria and the particle approximation. It should
be noted, however, that this is not a~difficulty in many applications
where the
probability of this event is zero. Nonetheless, to avoid this technical
problem, one natural strategy is to introduce randomized criteria
thresholds. We further assume that the parameters $(a_{n})_{n\geq0}$ are
sampled realizations of a collection of absolutely continuous random
variables $(A_{n})_{n\geq0}$. The main simplification of these randomized
criteria comes from the fact that the parameters $\epsilon_{m}:=\inf
_{0\leq
n\leq m}\inf_{t_{n}\leq s\leq t_{n+1}}|\mathcal
{H}_{t_{n},s}^{(n)}(\mathbb{P}%
_{\widehat{\eta}_{n},(t_{n},s)})-a_{n}|$\vspace*{1pt} are strictly positive for almost
every realization $A_{n}=a_{n}$ of the threshold parameters.\vspace*{-3pt}

\begin{theo}
For almost every realization of the random threshold parameters, and
for any
$\delta\in ( 0,1 ) $, we have the following exponential estimates:
\[
\mathbb{P} \bigl( \exists0\leq n\leq m\ t_{n}^{N}\not=t_{n}|(A_{n})_{0\leq
n\leq m} \bigr) \leq c_{1}(m) \biggl( 1+\delta\epsilon_{m}\sqrt{\frac{N}{2}%
} \biggr) \exp{ \bigl( - {N\delta^{2}\epsilon_{m}^{2}}/{%
c_{2}(m)} \bigr) }\vspace*{-3pt}\vadjust{\goodbreak}
\]
for some constants $c_{1}(m),c_{2}(m)<\infty$.
In addition, when the measures $H_{t_{n},s}^{(n)}$ have a finite support,
for any $\delta\in ( 0,1/(2\epsilon_{m}) ) $,
\[
\mathbb{P} \bigl( \exists0\leq n\leq mt_{n}^{N}\not=t_{n}|(A_{n})_{0\leq
n\leq m} \bigr) \leq c_{1}(m)\exp{ \bigl( - {N\delta
^{2}\epsilon_{m}^{2}}/{c_{2}(m)} \bigr) }
\]
holds for a possibly different pair of finite constants $%
c_{1}(m),c_{2}(m)<\infty$.
\end{theo}

\begin{pf}
Using Proposition~\ref{propreff}, we obtain the rather crude estimate
\begin{eqnarray*}
&&\mathbb{P} \bigl( \Omega-\Omega_{m}^{N} ( \delta,(A_{n})_{0\leq n\leq
m} ) |(A_{n})_{0\leq n\leq m}=(a_{n})_{0\leq n\leq m} \bigr) \\
&&\quad \leq\sum_{n=0}^{m}\sum_{s=t_{n}}^{t_{n+1}} \mathbb{P} \bigl( \bigl\vert
\mathcal{H}_{t_{n},s}^{(n)} \bigl( \mathbb{P}_{\widehat{\eta}%
_{n},(t_{n},s)}^{N} \bigr) -\mathcal{H}_{t_{n},s}^{(n)} \bigl( \mathbb{P}_{%
\widehat{\eta}_{n},(t_{n},s)} \bigr) \bigr\vert\geq\delta\bigl\vert
\mathcal{H}_{t_{n},s}^{(n)} \bigl( \mathbb{P}_{\widehat{\eta}%
_{n},(t_{n},s)} \bigr) -a_{n} \bigr\vert \bigr) \\
&&\quad \leq\sum_{n=0}^{m}\sum_{s=t_{n}}^{t_{n+1}}\mathbb{P} \bigl( \bigl\vert
\mathcal{H}_{t_{n},s}^{(n)} \bigl( \mathbb{P}_{\widehat{\eta}%
_{n},(t_{n},s)}^{N} \bigr) -\mathcal{H}_{t_{n},s}^{(n)} \bigl( \mathbb{P}_{%
\widehat{\eta}_{n},(t_{n},s)} \bigr) \bigr\vert
\geq\delta\epsilon_{m} \bigr) \\
&&\quad \leq c_{1}(m)\bigl(1+\delta\epsilon_{m}\sqrt{N/2}\bigr)\exp{ \bigl( -{N\delta
^{2}\epsilon_{m}^{2}}/{c_{2}(m)} \bigr) }%
\end{eqnarray*}
for a pair of finite constants $c_{1}(m),c_{2}(m)<\infty$. The final
line is
a direct consequence of Proposition~\ref{propreff} and an application of
Proposition~\ref{propke} completes the proof.
\end{pf}

We conclude that for almost every realization $(A_{n})_{0\leq n\leq
m}=(a_{n})_{0\leq n\leq m}$ the pair of particle models $(\mathcal{X}%
_{n}^{(N)},\widehat{\mathcal{X}}_{n}^{ ( N ) })_{0\leq n\leq m}$
and $(\mathcal{Y}^{(N)},\widehat{\mathcal{Y}}_{n}^{ ( N ) })_{0\leq
n\leq m}$ only differ on events $\Omega-\Omega_{m}^{N}(\delta
,(a_{n})_{0\leq
n\leq m})$ with exponentially small probabilities:
\begin{eqnarray*}
&&\mathbb{P} \bigl( \exists0\leq n\leq m\ \bigl(\mathcal{Y}^{(N)},\widehat{%
\mathcal{Y}}_{n}^{ ( N ) }\bigr)\not=\bigl(\mathcal{X}_{n}^{(N)},\widehat{%
\mathcal{X}}_{n}^{ ( N ) }\bigr)|(A_{n})_{0\leq n\leq
m}=(a_{n})_{0\leq n\leq m} \bigr)
\\
&&\quad \leq c_{1}(m)\bigl(1+\delta\epsilon_{m}\sqrt{N/2}\bigr)\exp{ \bigl( -{N\delta
^{2}\epsilon_{m}^{2}}/{c_{2}(m)} \bigr) }.
\end{eqnarray*}

\section{A functional central limit theorem}
\label{functtclsec}

\subsection{A direct approach}

In this section some direct consequences of the exponential coupling
estimates are discussed. For almost every realization $(A_{n})_{0\leq
n\leq
m}=(a_{n})_{0\leq n\leq m}$ and for any test function $f_{n}\in
\operatorname{Osc}_{1}(E_{t_{n}}^{\prime})$ the following
decomposition holds (writing $%
\overline{\eta}_{n}^{N}$ for the online adaptive approximation introduced
in Definition~\ref{defonline}):
\[
\sqrt{N} [ \overline{\eta}_{n}^{N}-\eta_{n} ] =\sqrt{N} [
\eta_{n}^{N}-\eta_{n} ] +\sqrt{N} [ \overline{\eta}%
_{n}^{N}-\eta_{n}^{N} ] 1_{\Omega-\Omega_{m}^{N}(\delta
,(a_{n})_{0\leq n\leq m})}
\]
with
\[
\mathbb{E} \bigl( \sqrt{N} [ \overline{\eta}_{n}^{N}-\eta_{n}^{N}%
] (f_{n})1_{\Omega-\Omega_{m}^{N}(\delta,(a_{n})_{0\leq n\leq
m})} \bigr) \leq\underbrace{\sqrt{N}\mathbb{P} ( \Omega-\Omega
_{m}^{N}(\delta,(a_{n})_{0\leq n\leq m}) ) } _{\stackrel{N\uparrow
\infty}{\longrightarrow}0}.
\]
Thus we can conclude directly that, for almost every realization $%
(A_{n})_{0\leq n\leq m}=(a_{n})_{0\leq n\leq m}$, the random fields
\[
\overline{W}_{n}^{N}:=\sqrt{N} [ \overline{\eta}_{n}^{N}-\eta_{n}%
] \quad \mbox{and}\quad  W_{n}^{N}:=\sqrt{N} [ \eta
_{n}^{N}-\eta_{n} ]\vspace*{-2pt}
\]
converge in law, as $N\uparrow\infty$, to the same centered Gaussian
random field $W_{n}$.

\vspace*{-2pt}\subsection{Functional central limit theorems}\vspace*{-2pt}

To demonstrate the impact of this functional fluctuation result we
provide a
brief discussion on the proof of the multivariate central limit
theorem. We
first recall the functional fluctuation theorem of the local errors
associated with the mean field particle approximation introduced in
(\ref%
{defVN}). This result was initially presented in~\cite{dmm} and
extended in~\cite{fk}.\vspace*{-3pt}

\begin{theo}
\label{keyth} For any fixed time horizon $n\geq0$, the sequence $%
(V_{p}^{N})_{0\leq p\leq n}$ converges in law, as $N$ tends to
infinity, to
a sequence of $n$ independent, Gaussian and centered random fields $%
(V_{p})_{0\leq p\leq n}$ with, for any $f_{p},g_{p}\in\mathcal{B}%
_{b}(E_{p}^{\prime})$, and $1\leq p\leq n$,
%
\begin{equation}\label{corr}
\mathbb{E}(V_{p}(f_{p})V_{p}(g_{p}))=\eta_{p-1}\mathcal{K}_{p,\eta
_{p-1}}\bigl([f_{p}-\mathcal{K}_{p,\eta_{p-1}}(f_{p})][g_{p}-\mathcal
{K}_{p,\eta
_{p-1}}(g_{p})]\bigr).\vspace*{-3pt}
\end{equation}
\end{theo}

Using arguments similar to those in the proof of Lemma~\ref{keydec}, we
obtain the decomposition formula:
\[
\lbrack\Phi_{n}(\mu)-\Phi_{n}(\eta)](f)=(\mu-\eta)\mathcal
{D}_{n,\eta
}(f)+\mathcal{R}_{n}(\mu,\eta)(f)\vspace*{-2pt}
\]
with the signed measure $\mathcal{R}_{n}(\mu,\eta)$ given by%
\begin{eqnarray*}
\mathcal{R}_{n}(\mu,\eta)(f) &:=&-\frac{1}{\mu(\mathcal
{G}_{n,\eta})}%
[\mu-\eta]^{\otimes2}\bigl(\mathcal{G}_{n,\eta}\otimes\mathcal
{D}_{n,\eta
}(f)\bigr)\qquad  \mbox{with } \mathcal{G}_{n,\eta}:=\mathcal{G}%
_{n-1}/\eta(\mathcal{G}_{n-1}), \\[-2pt]
\mathcal{D}_{n,\eta}(f)(x) &:=&\mathcal{G}_{n,\eta
}(x)\times%
\mathcal{M}_{n} \bigl( f-\Phi_{n}(\eta)(f) \bigr) (x).\vspace*{-3pt}
\end{eqnarray*}

\begin{defi}
Denote by $D_{p,n}$ the semi-group associated to the integral operators
$%
D_{n}:=\mathcal{D}_{n,\eta_{n-1}}$; that is, $D_{p,n}:=D_{p+1}\cdots
D_{n-1}D_{n}. $ For $p=n$, we use the convention $D_{n,n}=\mathit{Id}$, the identity
operator.\vspace*{-3pt}
\end{defi}

The semigroup $D_{p,n}$ can be explicitly described in terms of the
semigroup $\mathcal{Q }_{p,n}$ via
\[
D_{p,n}(f)=\frac{ \mathcal{Q }_{p,n}}{\eta_{p}( \mathcal{Q
}_{p,n}(1))}%
\bigl( f-\eta_{n}(f) \bigr) .\vspace*{-2pt}
\]

The next lemma provides a first-order decomposition of the random
fields $%
W_{n}^{N}$ in terms of the local fluctuation errors. Its proof is in the
\hyperref[app]{Appendix}. Note that $R_{p}$ can be understood in the proof.\vspace*{-2pt}

\begin{lem}
\label{lemdecomp} For any $N\geq1$ and any $0\leq p\leq n$, we have
%
\begin{equation} \label{RaN}
W^{N}_{n}=\sum_{p=0}^{n} V^{N}_{p} D_{p,n}+ \mathcal{R }^{N}_{n} %
\qquad \mbox{with } \mathcal{R }^{N}_{n}:=\sqrt{N} \sum_{p=0}^{n-1}
R_{p+1} ( \eta^{N}_{p},\eta_{p} ) D_{p+1,n}.\vspace*{-3pt}\vadjust{\goodbreak}
\end{equation}
\end{lem}

Using the $\mathbb{L}_{m}$-mean error estimates presented in
Section~\ref%
{lmbounds}, it is easily proved that the sequence of remainder random fields
$\mathcal{R }^{N}_{n}$ in (\ref{RaN}) converge in law, in the sense of
finite distributions, to the null random field as $N\uparrow\infty$.
Therefore the fluctuations of $W^{N}_{n}$ follow from Theorem~\ref{keyth}.

\begin{cor}
\label{lecor} For any fixed time horizon $n\geq0$, the sequence of random
fields $(W^{N}_{n})_{n\geq0}$ converges in law, as $N\uparrow\infty
$, to a
sequence of Gaussian and centered random fields $(W_{n})_{n\geq0}$,
where $%
\forall n\geq0$ $W_{n}=\sum_{p=0}^{n} V_{p}D_{p,n}. $
\end{cor}

\subsection{On the fluctuations of weighted occupation measures}

We end this article with some comments on the fluctuations of weighted
occupation measures on path spaces. Returning to the online adaptive
particle model, given $(t_{n}^{N},t_{n+1}^{N})=(t_{n},t_{n+1})$ the $N$%
-particle measures $\overline{\eta}_{n+1}^{N}=\frac{1}{N}%
\sum_{i=1}^{N}\delta_{ (\widehat{\mathcal{Y}}%
_{n}^{(N,i)},(Y_{t_{n}^{N}+1}^{(N,i)},Y_{t_{n}^{N}+2}^{(N,i)},\ldots
,Y_{t_{n+1}^{N}}^{(N,i)}) )}$ can be used to approximate the flow of
updated Feynman--Kac path distributions $(\widehat{\eta
}_{n+1,s})_{t_{n}\leq
s\leq t_{n+1}}$ given for any bounded test function $f_{n+1}\in
\mathcal{B}%
_{b}(S_{n+1})$ by
\[
s\in\lbrack t_{n},t_{n+1}]\mapsto\widehat{\eta
}_{n+1,s}(f_{n+1})\propto
\mathbb{E} [f_{n+1}(X_{0:t_{n+1}})W_{0:s} ( X_{1:s} )  ].
\]
Indeed, if we choose
\[
T_{n+1}^{(1)}(f_{n+1})(x_{0:t_{n+1}}):=f_{n+1}(x_{0:t_{n+1}})W_{t_{n}:s}
( x_{t_{n}+1:s} ),
\]
then in some sense
\[
\widehat{\eta}_{n+1,s}^{N}(f_{n+1}):=\frac{\overline{\eta
}_{n+1}^{N} (
T_{n+1}^{(1)}(f_{n+1}) ) }{\overline{\eta}_{n+1}^{N} (
T_{n+1}^{(1)}(1) ) }\simeq_{N\uparrow\infty}\widehat{\eta}%
_{n+1,s}(f_{n+1}):=\frac{\eta_{n+1} ( T_{n+1}^{(1)}(f_{n+1}) ) }{%
\eta_{n+1} ( T_{n+1}^{(1)}(1) ) },
\]
where $\eta_{n+1}$ is the flow of Feynman--Kac measures on path spaces
introduced in Section~\ref{secexc}.

Since the adaptive interaction time is taken such that
$t_{n+1}^{N}=t_{n+1}$%
, it holds that
\[
\frac{1}{N}\sum_{i=1}^{N}\delta_{\widehat{\mathcal
{Y}}_{n+1}^{(N,i)}}%
\simeq_{N\uparrow\infty}\widehat{\eta}_{n+1,t_{n+1}}=\widehat
{\eta}%
_{n+1}.
\]
In other words, if the marginal type functions are chosen such that
\[
T_{n+1}^{(0)}(f_{n+1})(x_{0:t_{n+2}}):=f_{n+1}(x_{0:t_{n+1}})
\]
so
\begin{eqnarray*}
\eta_{n+2} \bigl( T_{n+1}^{(0)}(f_{n+1}) \bigr) &=&\widehat{\eta}%
_{n+1}(f_{n+1})\propto\mathbb{E} [f_{n+1}(X_{0:t_{n+1}})W_{0:t_{n+1}}%
( X_{1:t_{n+1}} ) ], \\
\overline{\eta}_{n+2}^{N} \bigl( T_{n+1}^{(0)}(f_{n+1}) \bigr) &=&\frac{1}{N}%
\sum_{i=1}^{N}f_{n+1} \bigl( \widehat{\mathcal{Y}}_{n+1}^{(N,i)} \bigr)
\simeq_{N\uparrow\infty}\eta_{n+2} \bigl( T_{n+1}^{(0)}(f_{n+1}) \bigr) .
\end{eqnarray*}
From the previous discussion, for almost every realization
$(A_{n})_{0\leq
n\leq m}=(a_{n})_{0\leq n\leq m}$, a~central limit theorem (CLT) is
easily derived for the collection
of random fields
\begin{eqnarray*}
\widehat{W}_{n+1}^{N,(0)}(f_{n+1})&:= &\sqrt{N} \bigl[ \overline{\eta}%
_{n+2}^{N} \bigl( T_{n+1}^{(0)}(f_{n+1}) \bigr) -\eta_{n+2} \bigl(
T_{n+1}^{(0)}(f_{n+1}) \bigr) \bigr], \\
\widehat{W}_{n+1,s}^{N,(1)}(f_{n+1})&:= &\sqrt{N} [ \widehat{\eta}%
_{n+1,s}^{N}(f_{n+1})-\widehat{\eta}_{n+1,s}(f_{n+1}) ]
\end{eqnarray*}
as well as for the mixture of random field sequences
%
\begin{equation}
\widehat{W}_{n+1,s}^{N}=1_{t_{n}^{N}\leq s<t_{n+1}^{N}}\widehat{W}%
_{n+1,s}^{N,(1)}+1_{s=t_{n+1}^{N}}\widehat{W}_{n+1}^{N,(0)}.
\label{mixfield}
\end{equation}

The fluctuation analysis of these random fields relies on the
functional CLT
stated in Corollary~\ref{lecor}. In particular, the fluctuations of the
random fields (\ref{mixfield}) depend on those of a pair of random fields.

\subsection{Related work}

Reference \cite{douc} is the only published paper discussing a convergence
result for an adaptive SMC\ scheme. The authors establish a CLT using an
inductive proof w.r.t. deterministic time periods. They avoid the degenerate
situation where the threshold parameter coincides with the limiting
functional criterion. More recently, this problem has also been
addressed in
\cite{cornebise2009}, Chapter 4. However, the author does not account for
the randomness of the resampling times in his analysis.

\begin{appendix}
\section*{Appendix}\label{app}

\begin{pf*}{Proof of Lemma~\protect\ref{keydec}}
Via (\ref{phiQ}), for any $f\in\mathcal{B}_{b}(S_{n+1})$ we find that
\begin{eqnarray*}
\lbrack\Phi_{p,n}(\mu)-\Phi_{p,n}(\eta)](f)& \hspace*{3pt}= &\frac{1}{\mu
(\mathcal{G}%
_{p,n,\eta})}(\mu-\eta)\mathcal{D}_{p,n,\eta}(f), \\
\mathcal{D}_{p,n,\eta}(f)(x)& :=&\mathcal{G}_{p,n,\eta}(x)\times
\mathcal{P}%
_{p,n} \bigl( f-\Phi_{p,n}(\eta)(f) \bigr) (x),
\end{eqnarray*}
where $\mathcal{G}_{p,n,\eta}:=\mathcal{Q}_{p,n}(1)/\eta(\mathcal
{Q}%
_{p,n}(1))$ and $\mathcal{P}_{p,n}(f)=\mathcal{Q}_{p,n}(f)/\mathcal
{Q}%
_{p,n}(1)$. Now, since \mbox{$\eta( \mathcal{G}_{p,n,\eta} ) =1$}, it
follows that
\begin{eqnarray*}
\lbrack\Phi_{p,n}(\mu)-\Phi_{p,n}(\eta)]&\hspace*{3pt}=&(\mu-\eta)\mathcal
{D}_{p,n,\eta}+%
\mathcal{R}_{p,n}(\mu,\eta), \\
\mathcal{R}_{p,n}(\mu,\eta)(f)&:=&-\frac{1}{\mu(\mathcal
{G}_{p,n,\eta})}%
[\mu-\eta]^{\otimes2}\bigl(\mathcal{G}_{p,n,\eta}\otimes\mathcal
{D}_{p,n,\eta
}(f)\bigr).
\end{eqnarray*}
Using the fact that
\[
\mathcal{D}_{p,n,\eta}(f)(x)=\mathcal{G}_{p,n,\eta}(x)\int[
\mathcal{P%
}_{p,n}(f)(x)-\mathcal{P}_{p,n}(f)(y) ] \mathcal{G}_{p,n,\eta}(y)%
\eta(\mathrm{d}y)
\]
we find
\[
\forall f\in\operatorname{Osc}_{1}(S_{n})\qquad  \Vert\mathcal{D}%
_{p,n,\eta}(f)\Vert\leq q_{p,n}\beta(\mathcal{P}_{p,n}).
\]

Finally, for any $f\in\operatorname{Osc}_{1}(S_{n})$ observe that
\[
\vert\mathcal{R}_{p,n}(\mu,\eta)(f) \vert\leq (
4q_{p,n}^{3}\beta(\mathcal{P}_{p,n}) ) \bigl\vert[\mu-\eta
]^{\otimes2}\bigl ( \overline{\mathcal{G}}_{p,n,\eta}\otimes\overline{%
\mathcal{D}}_{p,n,\eta}(f) \bigr) \bigr\vert
\]
with $\overline{\mathcal{G}}_{p,n,\eta}:=\mathcal{G}_{p,n,\eta}/2q_{p,n}$
and $\overline{\mathcal{D}}_{p,n,\eta}(f):=D_{p,n,\eta
}(f)/2q_{p,n}\beta(%
\mathcal{P}_{p,n})\in\operatorname{Osc}_{1}(S_{p}).$
\end{pf*}

\begin{pf*}{Proof of Lemma~\protect\ref{lemdecomp}}
The lemma is proved by induction on $n$. For $n=0$, it follows that $%
W^{N}_{n}=V^{N}_{0}=\sqrt{N} [ \eta^{N}_{0}-\Phi_{0}(\eta^{N}_{-1})%
] $, with $\Phi_{0}(\eta^{N}_{-1})=\eta_{0}$. Assuming the formula
at $%
n$
\begin{eqnarray*}
W^{N}_{n+1} & = & V^{N}_{n+1}+\sqrt{N} [ \Phi_{n+1}(\eta^{N}_{n})-%
\Phi_{n+1}(\eta_{n}) ] \\
& = & V^{N}_{n+1}+W^{N}_{n}D_{n+1}+\sqrt{N}R_{n+1} (
\eta^{N}_{n},\eta_{n} ) \\
& = & V^{N}_{n+1}+\sum_{p=0}^{n} V^{N}_{p} D_{p,n+1} +\sqrt{N}%
\sum_{p=0}^{n-1} R_{p+1} ( \eta^{N}_{p},\eta_{p} ) D_{p+1,n+1}+%
\sqrt{N}R_{n+1} ( \eta^{N}_{n},\eta_{n} ) .
\end{eqnarray*}
Letting $D_{n+1,n+1}=I$, it follows that (\ref{RaN}) is satisfied at
rank $%
(n+1)$ due to
\begin{eqnarray*}
V^{N}_{n+1}+ \sum_{p=0}^{n} V^{N}_{p} D_{p,n+1} & = &%
\sum_{p=0}^{n+1} V^{N}_{p} D_{p,n+1}, \\
\sum_{p=0}^{n-1} R_{p+1} ( \eta^{N}_{p},\eta_{p} )
D_{p+1,n+1}+R_{n+1} ( \eta^{N}_{n},\eta_{n} ) & = & %
\sum_{p=0}^{n} R_{p+1} ( \eta^{N}_{p},\eta_{p} ) D_{p+1,n+1}.
\end{eqnarray*}
\upqed
\end{pf*}
\end{appendix}

\section*{Acknowledgements}

We would like to thank the Associate Editor and the three referees for many
helpful comments that vastly improved the paper.

%

\printhistory

\end{document}